\documentclass[12pt]{amsart}
\usepackage[latin1]{inputenc}
\usepackage{amsbsy}
\usepackage{amscd} 
\usepackage{amsfonts}
\usepackage{amsgen} 
\usepackage{amsmath}
\usepackage{amsopn} 
\usepackage{amssymb}
\usepackage{amstext}
\usepackage{amsthm} 
\usepackage{amsxtra}
\usepackage[mathscr]{eucal}

\theoremstyle{plain} 
\newtheorem{thm}{Theorem}[section]
\newtheorem{prop}[thm]{Proposition}
\newtheorem{lem}[thm]{Lemma}
\newtheorem{cor}[thm]{Corollary}
\theoremstyle{definition}
\newtheorem{defn}[thm]{Definition}
\newtheorem{rem}[thm]{Remark}

\newcommand{\Proof}{\textbf{Proof. }}

\numberwithin{equation}{section}

\renewcommand{\theta}{\vartheta}
\renewcommand{\phi}{\varphi}
\renewcommand{\epsilon}{\varepsilon}
\renewcommand{\subset}{\subseteq}
\renewcommand{\supset}{\supseteq}

\newcommand{\qqed}{\qed\medskip}
\newcommand{\N}{\mathbb N}
\newcommand{\Z}{\mathbb Z}
\newcommand{\Q}{\mathbb Q}
\newcommand{\R}{\mathbb R}
\newcommand{\C}{\mathbb C}
\newcommand{\F}{\mathbb F}
\newcommand{\T}{\mathbb T}
\newcommand{\K}{\mathcal K}
\newcommand{\norm}[1]{\lVert {#1}\rVert}
\newcommand{\scp}[2]{(#1 \, | \, #2)}
\newcommand{\surj}{\twoheadrightarrow}
\newcommand{\inj}{\hookrightarrow}
\newcommand{\eq}{\quad \Longleftrightarrow\quad}

\newcommand{\CS}{\mathcal C(S^1)}
\newcommand{\Toepl}{\mathscr{T}}
\newcommand{\ttens}[2]{#1 \otimes_\theta #2}
\newcommand{\Afrei}{\Toepl *_\theta \! \Toepl}
\newcommand{\Afreiminus}{\Toepl *_{-\theta} \! \Toepl}
\newcommand{\Afreimu}{\Toepl *_\mu \! \Toepl}
\newcommand{\Atens}{\ttens{\Toepl}{\Toepl}}
\newcommand{\Aver}{\Toepl \rtimes_\theta \Z}
\newcommand{\Averminus}{\Toepl \rtimes_{-\theta} \Z}
\newcommand{\Aup}{\CS *_\C \C^2}
\newcommand{\Asort}{(\Aver) \otimes (\Aup)}
\newcommand{\AsortT}{(\Aver)\otimes L}
\newcommand{\J}{(\Aup)\otimes \K}
\newcommand{\univ}[2]{C^*(#1 \, | \, #2)}
\newcommand{\univbig}[2]{C^*\!\big(#1 \, \big| \, #2\big)}
\newcommand{\LHH}[1]{\mathcal L(#1)}
\newcommand{\LN}{\mathcal L(\ell^2(\N_0))}
\newcommand{\LH}{\LHH{H}}
\newcommand{\On}[1]{\mathcal O_{#1}}
\newcommand{\starrow}[1]{\stackrel{#1}{\longrightarrow}}

\newcommand{\Iu}{\langle 1-uu^*\rangle}
\newcommand{\Iv}{\langle 1-vv^*\rangle}
\newcommand{\Iuv}{\langle 1-uu^*,1-vv^*\rangle}
\newcommand*{\lhdnicht}{%
\mathrel{\vcenter{\offinterlineskip
\vskip +.8ex\hbox{$\lhd$}\vskip -.2ex\hbox{$\neq$}}}}
\newcommand{\zw}[2]{\textbf{A.#1 #2.}}

\DeclareMathOperator{\id}{id}

\begin{document}
\title[Isometries with Twisted Commutation Relations]{On $C^*$-Algebras Generated by Isometries with Twisted Commutation Relations}
\author{Moritz Weber}
\address{Saarland University, Fachrichtung Mathematik, Postfach 151150,
66041 Saarbr\"ucken, Germany}
\email{weber@math.uni-sb.de}
\date{\today}
\subjclass[2010]{46L05 (Primary); 46L80, 19Kxx, 46L65 (Secondary)}
\keywords{universal $C^*$-algebra, rotation algebra, noncommutative torus, isometries, commutation relations, twist}
\thanks{\scriptsize{Research supported by the Deutsche Forschungsgemeinschaft and the Graduiertenkolleg ``Analytische Topologie und Metageometrie''. This work has been done in the context of the author's PhD project at the University of Muenster.}}

\begin{abstract}
In the theory of $C^*$-algebras, interesting noncommutative structures arise as deformations of the tensor product.
For instance, the rotation algebra $A_\theta$ 
may be seen as a scalar twist deformation of the tensor product $\CS\otimes\CS$.
We deform the tensor product 
of two Toeplitz algebras  in the same way, introducing 
the universal $C^*$-algebra $\Atens$ generated by two isometries $u$ and $v$ such that $uv=e^{2\pi i\theta}vu$  and $u^*v=e^{-2\pi i\theta}vu^*$, for $\theta\in\R$. Since the second relation implies the first one, we also consider 
the universal $C^*$-algebra $\Afrei$ generated by two isometries $u$ and $v$ with the weaker relation $uv=e^{2\pi i\theta}vu$. Such a "weaker case" does not exist in the case of unitaries, and it turns out to be much more interesting than the twisted "tensor product case" of two Toeplitz algebras.

We show that $\Atens$ is nuclear, whereas $\Afrei$ is not even exact. 
Also, we compute the $K$-groups and we obtain $K_0(\Afrei)=\Z$ and $K_1(\Afrei)=0$, and the same $K$-groups for $\Atens$. This answers a question raised by Murphy in 1994
concerning the $K$-theory of the $C^*$-algebra $C^*(\N^2)$, which is $\Afrei$ with $\theta=0$. 
\end{abstract}

\maketitle

\section*{Introduction}

\noindent
Given two unital $C^*$-algebras $A$ and $B$, a new object may be constructed by forming the (maximal) tensor product $A\otimes B$. It is given by the elements of $A$ and $B$ such that all elements of the first $C^*$-algebra commute with all elements of the second, i.e. $ab=ba$ for all $a\in A$ and $b\in B$. Deformations of these commutation relations (resp. of the tensor product construction) give rise to interesting structures in operator algebras, which are investigated with the aim of understanding settings with noncommutative multiplications.

An important example of a noncommutative $C^*$-algebra obtained in such a way is the (irrational) rotation algebra $A_\theta$. It has been studied since the early days of the theory (see for instance \cite{R81}, \cite{EH67}, \cite{Z68}), and it is one of the guiding lights in the development of several tools and theories (e.g. in Connes' noncommutative geometry, \cite{C85}, \cite{C94}). Furthermore, it serves as a model for many phenomena in physics (see \cite[Ch. 12]{GVF01} for a short overview).

The rotation algebra $A_\theta$ is the universal $C^*$-algebra generated by two unitaries $u$ and $v$ such that $uv=\lambda vu$, where $\lambda=e^{2\pi i \theta}$ is a scalar of absolute value one, for $\theta\in \R$. For $\theta=0$, the rotation algebra may be seen as the tensor product $\CS\otimes\CS$ of the algebra $\CS$ of continuous functions on the circle with itself. For arbitrary $\theta\in \R$ it is the twisted tensor product $\ttens{\CS}{\CS}$, i.e. it is the universal $C^*$-algebra generated by two unitaries with twisted commutation relations.

Unitaries are the isomorphisms in the category of Hilbert spaces -- they are isometric and surjective. Isometries (i.e. not necessarily surjective unitaries) in turn form a more general class of operators on  Hilbert spaces.
The universal $C^*$-algebra generated by a single isometry is the well-known Toeplitz algebra $\Toepl$, an extension of $\CS$ by the $C^*$-algebra $\K$ of compact operators.

In this article, we  investigate $C^*$-algebras generated by two isometries with twisted commutation relations. Surprisingly, a systematic study of this structure has not yet been done (for partial results by others, see the end of section \ref{SectTwist}). Furthermore, forming the $^*$-algebraic product of isometries reveals problems which do not arise in the case of unitaries.

A first step is to \emph{twist} the tensor product $\Toepl\otimes\Toepl$ of two Toeplitz algebras. This $C^*$-algebra is generated by two isometries $u$ and $v$ such that $uv=vu$ and $u^*v=vu^*$. 
The second relation implies the first one but not conversely (in contrast to the case of two unitaries). 
Therefore, we introduce two types of twisted commutation relations. First the \textbf{\emph{tensor twist of two isometries}} $\Atens$, i.e. the universal $C^*$-algebra generated by two isometries $u$ and $v$ such that $u^*v=\bar\lambda vu^*$ (and hence also $uv=\lambda vu$), where again $\lambda=e^{2\pi i \theta}$ and $\theta\in\R$. Secondly, the \textbf{\emph{free twist of two isometries}} $\Afrei$, i.e. the universal $C^*$-algebra generated by two isometries $u$ and $v$ such that $uv=\lambda vu$. 
For unitaries $u$ and $v$ there is only one way to twist the commutation relation $uv=vu$ whereas for isometries we have these two versions, which behave quite differently. It turns out that the $C^*$-algebra $\Afrei$ is in some sense much more complicated than $\Atens$.

Viewed on the complexity scale of nuclearity and exactness, $\Afrei$ carries much more ``free'' structure whereas $\Atens$ behaves more like the tensor product $\Toepl\otimes\Toepl$. This is due to the \textbf{first main result} of this article, namely that $\Atens$ is nuclear whereas $\Afrei$ is not even exact.  In order to show this second, much harder claim, we consider a Wold type decomposition of the isometries $u$ and $v$. We obtain a good description of the kernel of the canonical map $\Afrei\to A_\theta$ and of further important ideals in $\Afrei$. The decomposition is inspired by the work of Berger, Coburn and Lebow (\cite{BCL78}) on commuting isometries, which is $\Afrei$ in the case $\theta=0$.

This also leads to the \textbf{second main result}, the computation of the $K$-groups of $\Afrei$. It turns out that they are the same as those of $\Atens$. While the computations for the $C^*$-algebra $\Atens$ are not so hard, we have to make an effort in the case of $\Afrei$ going beyond tensor product like structures. The computation of the $K$-groups of $\Afrei$ answers a question raised by Murphy in 1994 (\cite{M94}) concerning the $C^*$-algebra $C^*(\N^2)$ which is $\Afrei$ in the case $\theta=0$.

Some aspects of $\Atens$ and $\Afrei$ were studied by Proskurin, Kabluchko, J\o rgensen, and Samo\u \i lenko (\cite{P00}, \cite{P00b}, \cite{K01}, \cite{JPS05}). See the end of section \ref{SectTwist} for an overview on their results and  other related work.\\

The article is organized as follows. In section \ref{SectTwist}, we introduce the $C^*$-algebras $\Atens$ and $\Afrei$. Section \ref{SectIdeale} deals with the ideals generated by a defect projection $1-uu^*$ resp. $1-vv^*$ in $\Atens$ (resp. $\Aver$). From this we deduce that $\Atens$ is nuclear. In section \ref{SectJ}, we prove that the ideal $J$ generated by the two defect projections $1-uu^*$ and $1-vv^*$ in $\Afrei$ contains all non-trivial ideals in $\Afrei$, for irrational parameters $\theta$. This yields the classification of the $C^*$-algebras $\Atens$ and $\Afrei$ with respect to the parameter $\theta$, if $\theta$ is irrational.

In section \ref{SectDecomp}, we construct an embedding $\iota$ of $\Afrei$ into $\Asort$ by analyzing the Wold decomposition of the isometry $uv\in\Afrei$. This enables us to give a satisfying picture of the ideal $J$ in $\Afrei$, namely it is isomorphic to $\J$. We infer that $\Afrei$ is not exact. For irrational parameters $\theta$, the ideal structure of $\Afrei$ corresponds to that of $\Aup$, which we investigate in some concrete cases, see section \ref{SectIdStr}. 

The $K$-theory of the $C^*$-algebras $\Atens$ and $\Afrei$ is computed in section \ref{SectKTheorie}. 
At the end of this article, an appendix on examples and constructions of universal $C^*$-algebras is attached.


\section{The twisted commutation relations of two isometries} \label{SectTwist}

\noindent
The rotation algebra $A_\theta$ is an important example of a noncommutative $C^*$-algebra (see the introduction of this article). Recall its definition as a universal $C$-algebra. 

\begin{defn}
 Let $\theta\in\R$ be a parameter and put $\lambda=e^{2\pi i \theta}$. The \emph{rotation algebra} $A_\theta$ is the universal $C^*$-algebra generated by two unitaries $u$ and $v$ such that $uv=\lambda vu$.
\end{defn}

One of the main features of the rotation algebra is its close relation to a commutative $C^*$-algebra, the algebra $\mathcal C(\T^2)$ of continuous functions on the 2-torus. This algebra may be seen as the universal $C^*$-algebra generated by two commuting unitaries, due to the Gelfand isomorphism.
For $\theta=0$, the rotation algebra and $\mathcal C(\mathbb T^2)$ coincide. Thus, the rotation algebra (also called the \emph{noncommutative torus}) is something of a soft step into the world of noncommutative $C^*$-algebras, and many questions about $A_\theta$ find the inspiration for their solutions in the case $\theta=0$. On the other hand, if the parameter $\theta$ is irrational, the rotation algebra is highly noncommutative.

Another perspective is to view the rotation algebra $A_\theta$ as the twisted tensor product of the algebra $\CS$ of continuous functions on the circle with itself. Indeed, as $\CS$ is the universal $C^*$-algebra generated by a single unitary, the tensor product $\CS\otimes\CS$ is the universal $C^*$-algebra generated by two commuting unitaries $u$ and $v$ (see the appendix for a short overview on universal $C^*$-algebras) which is $A_\theta$ in the case $\theta=0$. Hence, the rotation algebra $A_\theta$ is the deformed tensor product $\ttens{\CS}{\CS}$.

Our aim is to study versions of $A_\theta$ involving isometries instead of unitaries, which is a natural generalization. Isometries are isometric transformations of the underlying Hilbert space, i.e. they preserve the essential structure of the space, the inner product. Unitaries in turn are surjective isometries --  they are the isomorphisms in the category of Hilbert spaces. These properties may be expressed in an algebraic way, i.e. $u$ is a unitary if and only if $uu^*=u^*u=1$, whereas an isometry $v$ is given by $v^*v=1$.
The universal $C^*$-algebra generated by a single isometry $v$ is the Toeplitz algebra $\Toepl$, a well-known object. The canonical map from $\Toepl$ to $\CS$, mapping $v\mapsto u$, gives rise to  the following short exact sequence.
\[0\to\K\to\Toepl\to\CS\to 0\]
Here the $C^*$-algebra of compact operators on a separable Hilbert space are denoted by $\K$.
In order to twist the tensor product of two Toeplitz algebras, we observe that $\Toepl\otimes\Toepl$ is the universal $C^*$-algebra generated by two isometries $u$ and $v$ such that $uv=vu$ and $u^*v=vu^*$. (Note that the relation $u^*v=vu^*$ is needed to let all $^*$-monomials in $u$ commute with all $^*$-monomials in $v$.) If $u$ and $v$ are unitaries, these two relations are equivalent, but in the case of isometries we only have one of the implications. The following statement is also mentioned in an article by J\o rgensen, Proskurin, and Samo\u\i lenko (\cite{JPS05}).

\begin{lem} \label{USternV}
 Let $u$ and $v$ be two isometries in a unital $C^*$-algebra $A$ and let $\lambda \in S^1$ be a complex number of absolute value one. Then
\begin{enumerate}
 \item[(i)] $u^*v= \bar\lambda vu^* \Rightarrow uv=\lambda vu$
 \item[(ii)] $u^*v= \bar\lambda vu^* \not\Leftarrow uv=\lambda vu$ 
\end{enumerate}
\end{lem}
\Proof (i) We have $(uv - \lambda vu)^*(uv - \lambda vu) = 0$.

(ii) Consider $d(\lambda)\in \LN$, defined by $d(\lambda) e_n = \lambda^n e_n$, and the unilateral shift $S \in \LN$, defined by $S e_n = e_{n+1}$. Then $u':= d(\lambda)S$ and $v':= S$ are isometries with $u'v' = \lambda v' u'$, but $u'^*v' \not= \bar\lambda v' u'^*$ since $u'^*v' e_0 = \bar\lambda e_0$ whereas $v' u'^* e_0 = 0$.  
\qqed

Therefore, we distinguish two cases of twisted commutation relations.

\begin{defn}
Let  $\theta \in \R$ be a parameter and put $\lambda = e^{2 \pi i \theta} \in S^1$.
\begin{itemize}
 \item The \textit{tensor twist of two isometries}  $\Atens$ is defined as the universal $C^*$-algebra  generated  by two isometries $u$ and $v$ such that $u^*v = \bar\lambda vu^*$.
 \item  The \textit{free twist of two isometries} $\Afrei$ is defined as the universal $C^*$-algebra generated by two isometries $u$ and $v$ such that 
 $uv = \lambda vu$.
\end{itemize}
\end{defn}

We immediately see that $\Atens=\Toepl\otimes\Toepl$ if $\theta=0$. On the other hand, the (unital) free product of $\Toepl$ with itself and the $C^*$-algebra $\Afrei$ do \emph{not} coincide in the case $\theta=0$. Nevertheless, $\Afrei$ carries much more free structure than $\Atens$ since there are no direct relations between $u^*$ and $v$, i.e. the monomial structure of $\Afrei$ is much more complicated. Also we will see that the free group $C^*$-algebra $C^*(\F^2)$ embeds into $\Afrei$. This is a further hint about the features of a free product in the case of $\Afrei$.
In the following we will not distinguish between rational and irrational parameters $\theta$ if not explicitely stated. Also, we will always write $\lambda:= e^{2\pi i \theta}$.

Let us quickly remark that the deformed tensor product of $\CS\otimes\Toepl$ is a standard construction. Recall that $A_\theta$ may be constructed as a crossed product of $\CS$ by the group of integers $\Z$ according to the isomorphism $u\mapsto \lambda u$ on $\CS$. Likewise $\Aver$ is the crossed product of the Toeplitz algebra $\Toepl$ by the automorphism $v\mapsto \lambda v$.

\begin{rem}\label{DefAver}
 Let $\theta\in\R$ be a parameter and put $\lambda=e^{2\pi i \theta}$. Then $\Aver$ is the universal $C^*$-algebra generated by a unitary $u$ and an isometry $v$ such that $uv=\lambda vu$ (or equivalently $u^*v=\bar\lambda vu^*$).
\end{rem}

A first note on the difference between $\Atens$ and $\Afrei$ is the following remark on the range projections of the generating isometries.

\begin{rem} \label{ProjKomm}
Consider the range projections in $\Atens$ and $\Afrei$.
 \begin{itemize}
  \item[(a)] The range projections $uu^*$ and $vv^*$ in $\Atens$ commute.
  \item[(b)] The range projections $uu^*$ and $vv^*$ in $\Afrei$ do not commute.
 \end{itemize}
\end{rem}
\Proof A proof of (b) will be given later (see Remark \ref{ListeIota}).\qqed

At this point, we could ask whether the behavior of the commutator $[uu^*,vv^*]$ is the only difference between $\Atens$ and $\Afrei$. Let $B$ denote the quotient of $\Afrei$ by the ideal generated by the commutator $[uu^*,vv^*]$. By the representation in the proof of Lemma \ref{USternV}(b), $B$ does not equal $\Atens$, i.e. the canonical homomorphism from $\Atens$ to $B$ is not an isomorphism. Furthermore, we will see later that $\Atens$ is nuclear whereas $B$ is not even exact (see Proposition \ref{Kommutator}).\\

The $C^*$-algebras $\Atens$ and $\Afrei$ have been considered by several authors before, although from a different perspective and not in the homogeneous franework of twisting the commutation relations of isometries.
Proskurin (\cite{P00}, \cite{P00b}) studied the $C^*$-algebra $\Atens$ under the name $\mathcal A_{\{0\},\{\lambda_{ij}\}}$, and he proved that it is nuclear. We will give a different proof in section \ref{SectIdeale}. Furthermore he showed, that the defect ideal $\langle 1-uu^*,1-vv^*\rangle\lhd\Atens$ (which is the kernel of the map $\Atens\to A_\theta$) is the largest ideal in $\Atens$ for irrational parameter $\theta$. We extend his approach to the more general case of $\langle 1-uu^*,1-vv^*\rangle\lhd\Afrei$.
Proskurin also considered the Fock representation of $\Atens$ and he showed that it is faithful for irrational parameters. 

Kabluchko (\cite{K01}) extended Proskurin's result on the Fock representation to rational parameters. Using this, he gave a concrete description of the defect ideals of $\Atens$ (he denoted $\Atens$ by $\mathcal A_{\{0\},\Theta}$) for all parameters $\theta$, which can be found in section \ref{SectIdeale} -- although we prove it in a different way. 
He showed this using the Fock representation of $\Atens$; our proof is focused on the algebraic structures.

J\o rgensen, Proskurin and Samo\u\i lenko (\cite{JPS05}) proved in terms of a notion called ``$*$-wildness'', that $\Afrei$ -- which they call $\mathcal A_2^q$ -- is not nuclear. We even go further and show that $\Afrei$ is not exact. They also gave a classification of the $C^*$-algebras  $\Afrei$, depending on the parameter $\theta$. It is  directly transferred from the classification of the rotation algebras $A_\theta$. 

All of the above articles refer to an article by Bo\.zejko and Speicher (\cite{BS94}), where they introduced the so called $q_{ij}$-relations (or $q_{ij}$-CCR) on elements $d_1,\ldots,d_n$, namely
\[d_id_j^*-q_{ij}d_j^*d_i=\delta_{ij} \quad \textnormal{for } i,j=1,\ldots,n.\]
Here, the $q_{ij}$ are complex scalars with $|q_{ij}|\leq 1$ and $q_{ij}=\bar q_{ji}$.
For $n=2$, $q_{11}=q_{22}=0$ and $q_{12}=\bar\lambda$, these are the relations of $\Atens$ under the correspondence $u\leftrightarrow d_1^*$ and $v\leftrightarrow d_2^*$.
Thus, $\Atens$ may be seen as the universal $C^*$-algebra generated by $q_{ij}$-relations in a very special case. Note, that the $d_i$ are co-isometries in our case.

There is also some connection to Murphy's work \label{MurphysWork} on crossed products by semigroups. In his article \emph{Crossed products of $C^*$-algebras by semigroups of automorphisms} from 1994 (\cite{M94}), he defines (amongst other things) a $C^*$-algebra associated to (unital) semigroups equiped with a 2-cocycle. In \cite[Example 3.3]{M94}, he considers the $C^*$-algebra $C^*_\theta(\N^2)$, which is associated to the semigroup $\N^2$ and a 2-cocycle, constructed out of a single complex scalar $\lambda$ of absolute value one. This is the universal $C^*$-algebra generated by two isometries $u$ and $v$ fulfilling the relation $uv=\lambda vu$; hence $C^*_\theta(\N^2)=\Afrei$.

In the introduction to this article, he mentions that the $K$-theory of this object was unknown, even for the case of trivial $\lambda$, i.e. for $C^*(\N^2)=\Afrei$, where $\theta=0$. We may fill in this gap in chapter \ref{SectKTheorie} (see Remark \ref{MurphyKTheory}). According to Murphy, the knowledge of this $K$-theory would help in the theory of generalized Toeplitz operators (see \cite{M94} or \cite{M96} for references on this).

Murphy (\cite{M96}) also investigated the structure of the most important ideal in $\Afrei$, namely the defect ideal $\langle 1-uu^*,1-vv^*\rangle\lhd\Afrei$. He did this only for the case $\theta=0$, but not for arbitrary $\theta$.

Finally, we should mention the considerations by Berger, Coburn and Lebow (\cite{BCL78}) concerning the universal $C^*$-algebra generated by commuting isometries $u$ and $v$, thus the $C^*$-algebra $\Afrei$ in the case of $\theta=0$. Their work plays a crucial role in Section \ref{SectDecomp}. They analyzed the Wold decomposition of the isometry $uv\in\Afrei$ for $\theta=0$. We generalize it to arbitrary parameters $\theta$ which provides our main tool, namely a  transparent picture of $\Afrei$.

\section{Ideals in $\Atens$ and $\Aver$}\label{SectIdeale}

\noindent
In this section, we take a first look at the ideal structure of $\Atens$ and $\Aver$. We find a description of their defect ideals generated by $1-uu^*$ resp. by $1-vv^*$. 
The case of the $C$-algebra $\Afrei$ is more complicated and will be treated in the next two sections.
The main result of this section is to prove that $\Atens$ is nuclear.\\

From the universal property, we infer the existence of the following natural maps, mapping the generators $u\mapsto u$ and $v\mapsto v$.
\[\Atens\surj\Aver\surj A_\theta\]
The kernel of the map $\Aver\surj A_\theta$ is the ideal $\Iv$ generated by the defect projection of the isometry $v$ (recall Remark \ref{DefAver}). In $\Atens$, we have the ideals $\Iu$ and $\Iuv$ arising as kernels of the according maps to $\Aver$ resp. $A_\theta$. The ideal $\Iv$ in $\Atens$ is the kernel of the flipped quotient map $\Atens\to \Averminus$, given by $u\mapsto v$ and $v\mapsto u$.

In the case of $\Aver$ and $\Atens$, we may easily describe the defect ideals $\Iu$ and $\Iv$. They are given by a tensor product of $\CS$ resp. $\Toepl$ with the $C^*$-algebra of compact operators. To prove this, we need the following lemma on twisted tensor products with the compacts.
Note that the $C^*$-algebra $\K$ of compact operators may be seen as the universal $C^*$-algebra generated by elements $x_i$, $i\in \N_0$ with the relations $x_i^*x_j=\delta_{ij}$ for all $i,j\in \N_0$ (see also the appendix). The element $x_0$ is a minimal projection, whereas all $x_i$ are partial isometries.

\begin{defn} \label{DefTwistKomp}
Let $A$ be a $C^*$-algebra and $\alpha := (\alpha_i)_{i\in \N_0} \subset \textnormal{Aut}(A)$ be automorphisms $\alpha_i$ of $A$ for $i \in \N_0$ with $\alpha_0=\id$. Let $B_\alpha$ be the universal $C^*$-algebra generated by elements $a\in A$ (together with the relations of $A$) and elements $x_i$, $i\in \N_0$ such that $x_i^*x_j = \delta_{ij}x_0$ and $ax_i = x_i\alpha_i(a)$ for all $a\in A$, $i\in\N_0$. We define the \emph{twisted tensor product  $A \otimes_\alpha \K$ of $A$ with the compacts} to be the ideal generated by all products $ab$, $a \in A$, $b \in \K$ in $B_\alpha$. 
\end{defn}

If $A$ is unital, then $A \otimes_\alpha \K=B_\alpha$.
If $\alpha_i=\id$ for all $i\in\N_0$, then $A\otimes_\alpha\K$ coincides with the usual tensor product $A\otimes\K$ (see also the appendix). In fact, even for arbitrary automorphisms $\alpha_i$, the twisted tensor product $A\otimes_\alpha\K$ is isomorphic to the untwisted one, as  is shown in the next lemma.

\begin{lem} \label{TensKompUnTwist}
 Let $A$ be a $C^*$-algebra and $\alpha := (\alpha_i)_{i\in \N_0} \subset \textnormal{Aut}(A)$ be automorphisms on $A$, with $\alpha_0=\id$. The twisted tensor product $A \otimes_\alpha \K$ with the compacts is isomorphic to the untwisted tensor product:
\[A \otimes_\alpha \K \cong A \otimes \K\qquad \textnormal{via } \; ax_i \leftrightarrow \alpha_i(a) \otimes x_i\]
\end{lem}
\Proof We write $A \otimes \K$ as $A \otimes_{\id} \K \lhd B_{\id}$ with the automorphisms $\alpha_i = \id$.
 Let $\pi: B_{\id} \to \LH$ be a faithful representation of $B_{\id}$ and decompose $H$ by $H = \left(\bigoplus_{i \geq 0} \pi (x_i x_i^*)H\right) \oplus K$ for a subspace $K \subset H$ (recall that the $x_i x_i^*$ are mutually orthogonal projections). 
 As $A$ may be embedded into $B_{\id}$ and since $ax_ix_i^*=x_ix_i^*a$ for all $a\in A$, $i\in\N_0$, we get representations $\tilde\beta_i : A \to \LHH{\pi(x_ix_i^*)H}$, $a \mapsto \pi(\alpha_i(a))\pi(x_ix_i^*)$ for $i\in \N_0$; likewise
$\tilde \beta_K: A \to \LHH{K}$, $a\mapsto \pi(a)p_K$,  
where $p_K$ is the projection onto $K\subset H$.  Thus, we get a representation 
$\tilde \beta:A \to \LH$, $a\mapsto \bigoplus_{i\geq 0} \tilde \beta_i(a) \oplus \tilde \beta_K(a)$ of $A$ on $\LH$. It fulfills 
$\tilde\beta(a)\pi(x_i)=\pi(\alpha_i(a)x_i)\in\pi(A\otimes_{\id}\K)$ and 
 $\pi(x_i)\tilde\beta(\alpha_i(a))=\pi(x_ix_0)\tilde\beta_0(\alpha_i(a)x_0)=\pi(\alpha_i(a)x_i)=\tilde \beta(a)\pi(x_i)$.

From this, we get a representation $\sigma:B_\alpha \to \LH$, mapping $a\mapsto \tilde\beta(a)$ and $x_i\mapsto \pi(x_i)$. Restricting $\sigma$ to $A\otimes_\alpha\K$, we obtain a map from $A\otimes_\alpha\K$ to $A \otimes_{\id}\K$.
Note that $A\otimes_\alpha\K$ is spanned by elements $ax_ix_j^*$, hence 
 the image of $A\otimes_\alpha\K$ under $\sigma$ is in $\pi(A\otimes_{\id}\K)$ since
$\sigma(ax_ix_j^*)=\pi(\alpha_i(a)x_i)\pi(x_j^*)\in\pi(A\otimes_{\id}\K)$.

So, we constructed a homomorphism $A\otimes_\alpha\K \to A\otimes_{\id} \K, ax_i\mapsto\alpha_i(a)x_i$, and likewise we obtain a homomorphism 
$A\otimes_{\id}\K \to A\otimes_\alpha \K, ax_i\mapsto\alpha_i^{-1}(a)x_i$.
These maps  are inverse to each other. \qqed

Recall that $\CS$ is the universal $C^*$-algebra generated by a unitary $u$, whereas $\Toepl$ is the universal $C^*$-algebra generated by an isometry $v$.

\begin{prop} \label{DefIdealAverAtens}
 The ideals generated by the single defect projections in $\Aver$ and $\Atens$ are of the following form.  
 \begin{itemize}
  \item[(a)] $\CS \otimes \K \cong \langle 1-vv^* \rangle \lhd \Aver$ \, via \, $u\otimes x_i \mapsto v^iu(1-vv^*)$
  \item[(b)] $\Toepl \otimes \K \cong \langle 1-vv^* \rangle \lhd \Atens$  \, via \, $v\otimes x_i \mapsto v^iu(1-vv^*)$
  \item[] $\Toepl \otimes \K\cong \langle 1-uu^* \rangle \lhd \Atens$   \, via \,  $v\otimes x_i \mapsto u^iv(1-uu^*)$
 \end{itemize}
\end{prop}
\Proof (a) We consider the automorphisms $\alpha_j$ on $\CS$ for $j\in\N_0$, given by  $u\mapsto \lambda^ju$, where $\lambda=e^{2\pi i \theta}$. Hence, we can form the according twisted tensor product $\CS\otimes_\theta\K$ (as an ideal in $B_\theta$) which is isomorphic to $\CS\otimes\K$ by the preceding lemma.  Here, $B_\theta$ is the universal $C^*$-algebra generated by elements $x_i$, $i\in \N_0$ and a unitary $w$ such that $x_i^*x_j=\delta_{ij}x_0$ and $wx_i=\lambda^ix_iw$ for all $i\in\N_0$.

The elements $x_i':=v^i(1-vv^*)$ and $u$ in $\Aver$ fulfill the relations $x_i'^*x_j'=\delta_{ij}x_0'$ and $ux_i'=\lambda^ix_i'u$ for all $i\in\N_0$. Thus, there is a homomorphism $\beta:B_\theta\to\Aver$ mapping $x_i\mapsto x_i'$ and $w\mapsto u$. 
The ideal $\ttens{\CS}{\K}$  in $B_\theta$ is spanned by all elements $w^kx_ix_j^*$ for $k\in\Z$, $i,j\in\N_0$. Therefore, the restriction of $\beta$ to $\ttens{\CS}{\K}$ maps onto the closed linear span $\mathcal S$ of all elements $\beta(w^kx_ix_j^*)=u^kv^i(1-vv^*)(v^*)^j$ for $k\in\Z$, $i,j\in\N_0$. The span $\mathcal S$ is an ideal in $\Aver$ containing the projection $1-vv^*$, hence $\beta(\ttens{\CS}{\K})=\mathcal S=\Iv$.

We next show that the restriction of $\beta$ to $\ttens{\CS}{\K}$ is injective. Choose a unital, faithful representation $\pi:B_\theta\inj\LH$ and decompose $H=\left(\bigoplus_{i\geq0}\pi(x_ix_i^*)H\right)\oplus K$, where $K\subset H$ is a subspace of $H$. 
Write $(e_i)_{i\geq0}$ for the standard orthonormal basis of $\ell^2(\N_0)$. There is a Hilbert space isomorphism $\bigoplus_{i\geq 0} \pi(x_ix_i^*)H \cong \ell^2(\N_0)\otimes H_0$, with $H_0 := \pi(x_0)H$, via the mapping $\pi(x_i)\eta \leftrightarrow e_i\otimes\eta$. (Note, that $\pi(x_i^*)H=H_0$, because $\pi(x_i^*)=\pi((x_ix_0)^*)=\pi(x_0)\pi(x_i^*)$ and $\pi(x_0)=\pi(x_i^*x_i)$.) 
Now, we analyze $\pi(w)$ and $\pi(x_i)$ under this decomposition of $H$.

The subspace $H_0\subset H$ is invariant under $\pi(w)$, since $\pi(w)\eta=\pi(wx_0)\eta=\pi(x_0w)\eta$ is in $H_0$ for all $\eta \in H_0$. Even more, $\pi(w)H_0=H_0$ holds, since $\pi(w^*)H_0\subset H_0$. Therefore $\tilde w:= \pi(w)_{|H_0} \in \LHH{H_0}$ is a unitary.

Under the Hilbert space isomorphism of $K^\perp$ and $\ell^2(\N_0)\otimes H_0$, the operator  $\pi(w)_{|K^\perp}$ corresponds to $d(\lambda)\otimes \tilde w$, where $d(\lambda)\in\LN$ denotes the diagonal operator given by $d(\lambda)e_n=\lambda^n e_n$. This is due to the following computation:
\begin{gather*}
 \pi(w)_{|K^\perp}(e_i\otimes\eta) \leftrightarrow \pi(w)\pi(x_i)\eta=\lambda^i\pi(x_i)\pi(w)\eta \leftrightarrow (d(\lambda)\otimes\tilde w)(e_i\otimes\eta)
\end{gather*}
The operator $\pi(x_j)_{|K^\perp}$ correponds to $S^j(1-SS^*)\otimes 1$, where $S\in\LN$ is the unilateral shift, given by $Se_n=e_{n+1}$:
\begin{align*}
 \pi(x_j)_{|K^\perp}(e_i\otimes\eta) &\leftrightarrow \pi(x_jx_i)\eta=\pi(x_jx_0^*x_i)\eta=\delta_{i0}\pi(x_j)\eta \\
&\leftrightarrow (S^j(1-SS^*)\otimes 1)(e_i\otimes\eta)
\end{align*}
On the subspace $K\subset H$, $\pi(w)_{|K}$ is a unitary, since $\pi(w)K=K$.  Furthermore, the operators $\pi(x_i)$ get annihilated on $K$ for all $i\in\N_0$. Summarizing, we can write 
\[\pi(w)=(d(\lambda)\otimes \tilde w) \oplus \pi(w)_{|K} \qquad \textnormal{and} \qquad \pi(x_i)=(S^i(1-SS^*)\otimes 1)\oplus 0.\] 

The restriction of $\pi$ to $\ttens{\CS}{\K}$ is a  representation on the Hilbert space $K^\perp$, because $\pi(w^kx_ix_j^*)=(d(\lambda)^kS^i(1-SS^*){S^*}^j\otimes\tilde w^k)\oplus 0$. On the other hand, there is a map $\sigma: \Aver \to \LHH{K^\perp}$, mapping $u\mapsto d(\lambda)\otimes \tilde w$ and $v \mapsto S\otimes 1$. Hence, the following diagram commutes and the restriction of $\beta$ to $\ttens{\CS}{\K}$ is injective (since $\pi$ is injective).
\begin{align*}
 \ttens{\CS}{\K} \qquad  &\stackrel{\beta}{\to}\qquad \Aver \\
\stackrel{\searrow}{\pi} \qquad &\LHH{K^\perp} \quad \stackrel{\swarrow}{\sigma} 
\end{align*}

We conclude that $\ttens{\CS}{\K}$ is isomorphic to the ideal $\Iv$ in $\Aver$. The ideal $\ttens{\CS}{\K}$ in $B_\theta$ in turn is isomorphic to $\CS\otimes\K$ via $wx_i\mapsto \lambda^iu\otimes x_i$ by Lemma \ref{TensKompUnTwist}.

(b) The proof for the $C^*$-algebra $\Atens$ is exactly the same. Note that in this case, the operators $\tilde w \in \LHH{H_0}$ and $\pi(w)_{|K} \in \LHH{K}$ are not unitaries but isometries.
\qqed

A proof of the fact that $\Toepl \otimes \K$ is isomorphic to the defect ideal 
$\langle 1-vv^* \rangle \lhd \Atens$ may be found in an article by Kabluchko (\cite[prop. 4]{K01}), but he proves it in a quite different way by means of representation theory, via a Fock representation. 
We attack the problem from a different perspective analyzing the algebraic structure of the span and  ``untwisting'' the twist with the compacts.

The result of Proposition \ref{DefIdealAverAtens}(i) may also be deduced from a result by  Williams (\cite{W07}). 
We form the crossed product of the following short exact sequence with the automorphism $v\mapsto \lambda v$ on $\Toepl$.
\[ 0\longrightarrow \K \longrightarrow \Toepl \longrightarrow \CS\longrightarrow 0 \]
Hence, we obtain the following exact sequence.
\[ 0\longrightarrow \K\rtimes_\theta \Z \longrightarrow \Aver \longrightarrow \CS \rtimes_\theta \Z\longrightarrow 0 \]
The kernel of the map $\Aver \to \CS \rtimes_\theta \Z=A_\theta$ is the defect ideal $\langle 1-vv^* \rangle$, on the other hand $\K\rtimes_\theta\Z \cong \CS\otimes\K$. 

As an immediate consequence of the preceding proposition, we obtain the following result on $\Atens$.

\begin{thm} \label{AtensNucl}
The $C^*$-algebra $\Atens$ is nuclear.
\end{thm}
\Proof By the previous proposition we have the following short exact sequence:
\[0\to\Toepl\otimes\K\to\Atens\to\Aver\to 0\]
Since nuclearity is preserved under taking crossed products with amenable groups, we infer that $\Aver$ is nuclear. Now, $\Toepl\otimes\K$ is nuclear, too, and hence also $\Atens$.\qqed

We  now sketch a more direct proof of the fact that $\Atens$ is nuclear using a result by Rosenberg from 1977 (\cite{R77}). 

\begin{lem}[{\cite[th. 3]{R77}}]\label{LemRosenberg}
 Let $A$ be a unital $C^*$-algebra, $B\subset A$ be a nuclear $C^*$-subalgebra containing the unit of $A$, $s\in A$ be an isometry such that $sBs^*\subset B$ and let $A=C^*(B,s)$, i.e. let $A$ be generated by $B$ and $s$. Then $A$ is nuclear.
\end{lem} 

In fact, the isometry $s\in A$ induces an endomorphism $\textnormal{Ad}(s):B\to B$, since $sBs^*\subset B$. Thus, Rosenberg's lemma states that nuclearity is  preserved under taking crossed products with the semigroup $\N$, since $A=B\rtimes_{\textnormal{Ad}(s)}\N$.

Now, the $C^*$-subalgebra $B_1:=\univ{u^kv^n(v^*)^n(u^*)^k}{n,k\in\N_0} \subset \Atens$ is commutative and hence nuclear. We apply Rosenberg's lemma to $B_1$, the isometry $u$, and $A_1:=C^*(B_1,u)\subset\Atens$. Thus, $A_1$ is nuclear. We apply it again to $B_2:=A_1$ and $v\in C^*(B_2,v)=\Atens$, and we infer that $\Atens$ is nuclear. This proof is a slight modification of a proof by Proskurin (\cite[Prop. 3]{P00}, \cite{P00b}).

A glance at the details of the sketched proof shows that the quotient $D$ of $\Afrei$ by the ideal generated by all commutators $[u^a(u^*)^a, v^b(v^*)^b]$, $a,b\in \N_0$ is nuclear. This is because the $C^*$-subalgebra $B_1\subset D$ -- defined analogously to  $B_1\subset\Atens$ -- is commutative, exactly because of the relations $[u^a(u^*)^a, v^b(v^*)^b]=0$ for all $a, b\in\N_0$.
Nevertheless, the canonical map from $D$ to $\Atens$ is not an isomorphism, again due to the representation of Lemma \ref{USternV}(b). This explains again that the commutators $[u^a(u^*)^a, v^b(v^*)^b]$ are not the only difference between $\Atens$ and $\Afrei$, refining the answer to the question after Remark \ref{ProjKomm}.

The $C^*$-algebra $\Afrei$ in turn is \emph{not} nuclear.  This was shown by J\o rgensen, Proskurin and Samo\u\i lenko within their concept of  ``$^*$-wildness''(\cite{JPS05}). A $C^*$-algebra is called $^*$-wild, if its representation theory may be traced back in a certain way (which we do not specify here) to the representation theory of the $C^*$-algebra $C^*(\F_2)$ of the free group with two generators. They showed that $\Afrei$ is $^*$-wild and that every $^*$-wild $C^*$-algebra is not nuclear. 
In section \ref{SectDecomp} we will show that $\Afrei$ is not even exact.
A satisifying study of the ideals in $\Afrei$ and also of the ideal $\Iuv$ in $\Atens$ require further preparation. An explicit description is given in section \ref{SectDecomp}.

We end this section by a remark on the ideal $\langle (1-uu^*)(1-vv^*)\rangle$ in $\Atens$ and on the quotient by it. In some sense, it shows how close $\Atens$ and $\Aver$ are. Again, we view $\Aver$ as the universal $C^*$-algebra generated by a unitary $u$ and an isometry $v$ such that $uv=\lambda vu$ (see Remark \ref{DefAver}).

\begin{prop} \label{AtensUndAver}
 The ideal $\langle (1-uu^*)(1-vv^*)\rangle$ in $\Atens$ is isomorphic to the $C^*$-algebra $\K$ of compact operators. The quotient $D$ of $\Atens$ by this ideal embeds into $(\Aver)\oplus(\Aver)$ via $u\mapsto u\oplus v$ and $v\mapsto v\oplus u^*$. 
\end{prop}
\Proof The ideal $\langle (1-uu^*)(1-vv^*)\rangle$ is the closed linear span of all elements of the form $e_{(a,b)(d,c)}:=u^av^b(1-uu^*)(1-vv^*)(v^*)^c(u^*)^d$ for $a,b,c,d\in\N_0$. These elements fulfill the relations of the matrix units, i.e. $e_{(a,b)(d,c)}^*=e_{(d,c)(a,b)}$ and $e_{(a,b)(d,c)}e_{(f,g)(i,h)}=\delta_{(d,c)(f,g)}e_{(a,b)(i,h)}$, where $\delta_{(d,c)(f,g)}$ is defined as the Kronecker delta $\delta_{df}\delta_{cg}$. Thus $\langle (1-uu^*)(1-vv^*)\rangle$ is isomorphic to $\K$ (see also the appendix).

By the universal property, there is a homomorphism $\phi$ from the quotient $D$ to $(\Aver)\oplus(\Aver)$ such that $u\mapsto u\oplus v$ and $v\mapsto v\oplus u^*$. Furthermore, let $\pi$ be a unital, faithful representation of $D$ on a Hilbert space $H$. Then $H$ may be decomposed into $H=H_1\oplus H_2\oplus H_3$ such that 
\begin{itemize}
 \item $\pi(u)=u_1\oplus u_2\oplus u_3$, where $u_1$ is an isometry and $u_2$, $u_3$ are unitaries,
 \item $\pi(v)=v_1\oplus v_2\oplus v_3$, where $v_1$, $v_2$ are unitaries and $v_3$ is an isometry,
 \item and $u_iv_i=\lambda v_iu_i$ for $i=1,2,3$.
\end{itemize}
Indeed, put  $a_i:=(uv)^i(1-uu^*)((uv)^*)^i$ and $b_j:=(uv)^j(1-vv^*)((uv)^*)^j$ in $D$ for $i,j\in\N_0$. Then the $a_i$ and $b_j$ are mutually orthogonal projections.
We also have $ua_i=a_{i+1}u$, $ub_i=b_iu$, $va_i=a_iv$ and $vb_i=b_{i+1}v$ for all $i\in\N_0$. Check for instance:
\begin{align*}
ua_i&=\lambda^i(uv)^iu(vv^*+(1-vv^*))(1-uu^*)((uv)^*)^iu^*u \\
 &=\lambda^i (uv)^{i+1}v^*(1-uu^*)((uv)^*)^iu^*u \\
 &=\lambda^i \lambda^{-i} (uv)^{i+1}(1-uu^*)((uv)^*)^{i+1}u \\
 &=a_{i+1}u
\end{align*}

We put $H_1:=\bigoplus_{i\in\N_0}\pi(a_i) H$, $H_3:=\bigoplus_{i\in\N_0}\pi(b_i) H$ and $H_2:= H\ominus (H_1\oplus H_3)$. We infer from the relations on $u$, $a_i$ and $b_i$ that $\pi(u)H_1\subset H_1$,  $\pi(u)H_2=H_2$, and $\pi(u)H_3=H_3$. Similarly we get $\pi(v)H_1=H_1$, $\pi(v)H_2=H_2$ and $\pi(v)H_3\subset H_3$. 

We have a homomorphism $\sigma:(\Aver)\oplus(\Aver)\to\LH$, mapping
\begin{align*}
 &u\oplus 0\mapsto 0\oplus u_2\oplus u_3 && 0\oplus u\mapsto v_1^*\oplus0\oplus0\\
 &v\oplus 0\mapsto 0\oplus v_2\oplus v_3 && 0\oplus v\mapsto u_1\oplus0\oplus0
\end{align*}
Since $\sigma\circ\phi=\pi$, we conclude that $\phi$ is injective.
\qed


\section{The kernel of the map from $\Afrei$ to $A_\theta$} \label{SectJ}

\noindent
The kernel of the canonical map from $\Afrei$ to $A_\theta$ mapping $u\mapsto u$ and $v\mapsto v$ is the ideal $\Iuv$ generated by the defect projections $1-uu^*$ and $1-vv^*$. It plays an important role in the sequel and it is denoted by $J$. We prove that it reflects the whole ideal structure of $\Afrei$, if $\theta$ is irrational. Also, we give a classification of $\Atens$ and $\Afrei$ with respect to the parameter $\theta$.\\

Since $A_\theta$ is simple for irrational parameters $\theta$, we know that $J$ is a maximal ideal in $\Afrei$. Even more, we can show that it contains all non-trivial ideals in $\Afrei$. We need a technical lemma to prove this statement.

\begin{lem} \label{H2}
 Let $\theta$ be irrational, $I\lhd\Afrei$ be an ideal in $\Afrei$, $\epsilon>0$ and let $1=w+y+z$ be a decomposition of the unit such that
 \begin{itemize}
   \item $w$ is in the linear span of all elements $x(1-vv^*)x'$, where $x$ and $x'$ are $^*$-monomials in $u$ and $v$,
   \item $y\in I$,
   \item $z\in\Afrei$ with $\|z\|<\epsilon$.
 \end{itemize}
Then there exists a $y'\in I$, such that $\|y'-1\| < \epsilon$.

We obtain a similar result, if $w$ is in the linear span of all elements $x(1-uu^*)x'$. In fact, this span is dense in the ideal $\Iu\lhd\Afrei$.
\end{lem}
\Proof Put $s:=u^Nv^N$, where $N\in\N$ is sufficiently large such that $s^*w=0$. Put $y':=s^*ys$. Then $1=s^*(w+y+z)s=y'+s^*zs$ and thus $\|y'-1\|\leq\|z\|<\epsilon$.
\qqed

\begin{prop} \label{JGroesst}
 If $\theta$ is irrational, then the ideal $J:=\langle 1-uu^*,1-vv^*\rangle$ contains all non-trivial ideals in $\Afrei$, i.e. for any ideal $\Afrei \not = I\lhd\Afrei$ we have $I\subset J$.  
\end{prop}
\Proof Let $I\lhd\Afrei$ be a non-zero ideal. Since $A_\theta$ is simple, we either have $I\subset J$ or $I+J=\Afrei$. In the latter case, we have $1\in I+J=I'+\Iv$ for $I':=I+\Iu\lhd\Afrei$. Hence, there is a $x\in \Iv$ and a $y\in I'$ such that $1=x+y$. For $m\in\N$ and $\epsilon_m:=\frac{1}{m}$, there are $w_m$ and $z_m$ in $\Afrei$ such that $1=w_m+y+z_m$ is a decomposition of the unit in the sense of the preceding lemma, where $x=w_m + z_m$. 
We obtain a sequence of elements $y_m'\in I'$ such that $y_m'\to 1$ for $m \to \infty$. Applying the lemma again to $1\in I'=I+\Iu$ yields a sequence $y_m''\in I$ such that $y_m''\to 1$ for $m\to\infty$. We conclude that $I=\Afrei$ whenever $I\not\subset J$. \qqed

The proof of Proposition \ref{JGroesst} is adapted from the proof for a simpler case, namely for $\Atens$ given by Proskurin (\cite[prop. 6]{P00}, \cite{P00b}). In $\Atens$, the ideal $\langle 1-uu^*, 1-vv^*\rangle$ has a simple form: it is the closed linear span of elements of the form:
\[u^av^b(1-uu^*)^{\epsilon_1}(1-vv^*)^{\epsilon_2}(v^*)^c(u^*)^d \]
Here $a,b,c,d\in\N_0, \epsilon_1,\epsilon_2\in\{0,1\}$ and $\epsilon_1+\epsilon_2\not = 0$.
In $\Afrei$ however, the ideal $\Iuv$ is more complicated.

As a corollary from Proposition \ref{JGroesst}, we get that the ideal generated by the two defect projections is the largest ideal in \emph{any} $C^*$-algebra, which is generated by two  isometries $u$ and $v$ with $uv=\lambda vu$, whenever $\theta\in\R\backslash\Q$.

\begin{cor}\label{allgJGroesst}
 Let $A$ be a unital, non-zero $C^*$-algebra and let $u,v\in A$ be two isometries with $uv=\lambda vu$, where $\lambda=e^{2\pi i\theta}$, $\theta\in\R\backslash\Q$. 

Then $\langle 1-uu^*, 1-vv^*\rangle\lhd C^*(u,v)\subset A$ is the union of all non-trivial ideals in the $C^*$-subalgebra generated by $u$ and $v$.
\end{cor}
\Proof Let $\phi:\Afrei\surj C^*(u,v)=:B\subset A$ be the homomorphism with $u\mapsto u, v\mapsto v$. Then $\phi(J)=J':=\langle 1-uu^*,1-vv^*\rangle\lhd B$. Let $B\not = I'\lhd B$ be an ideal in $B$. Its preimage $I:=\phi^{-1}(I')\lhd\Afrei$ is an  ideal in $\Afrei$, thus $I\subset J$ and hence $I'=\phi(I)\subset \phi(J)=J'$. \qqed

\begin{cor} \label{JGroesstInAtens}
Let $\theta$ be irrational.
 The ideal $J'=\langle 1-uu^*, 1-vv^*\rangle$ is the union of all non-trivial ideals in $\Atens$. The same holds true for $\Iv$ in $\Aver$.
\end{cor}

From the fact that $J\lhd\Afrei$ is the largest ideal in $\Afrei$, we may infer their classification with respect to the parameter $\theta$ (for irrational $\theta$). The following proof is mainly taken from \cite{JPS05}.

\begin{prop}[{\cite[prop. 2]{JPS05}}]
 For irrational parameter $\theta$, there is a unique, normalized trace on $\Afrei$ as well as on $\Atens$, induced by the trace on $A_\theta$. 
\end{prop}
\Proof On the irrational rotation algebra $A_\theta$ there is a unique, normalized trace $\tau:A_\theta \to \C$. 
Thus, $\tilde \tau:= \tau \circ \phi$ is a normalized trace on $\Afrei$, where $\phi:\Afrei\to A_\theta$ is the canonical homomorphism.

Let $\tilde \sigma$ be another normalized trace on $\Afrei$. Then $\tilde\sigma(J)=0$, because 
\[ |\tilde\sigma(x(1-uu^*)x')|^2 = |\tilde\sigma((1-uu^*)x'x)|^2\leq\tilde\sigma(1-uu^*)\tilde\sigma(x^*x'^*x'x)=0\]
for all monomials $x, x'$ in $\Afrei$ by Cauchy-Schwarz; likewise $|\tilde\sigma(x(1-vv^*)x')|^2=0$.
Then $\sigma:A_\theta\to\C$ is a trace on $A_\theta$, given by $\sigma(b):=\tilde\sigma(a)$ for an $a\in \Afrei$ such that $\phi(a)=b$. Therefore we have $\sigma=\tau$, which implies $\tilde\sigma=\sigma\circ\phi=\tau\circ\phi=\tilde\tau$.\qqed

J\o rgensen, Proskurin and Samo\u\i lenko concluded in their article \cite{JPS05} that the classification of the $C^*$-algebras $\Afrei$ is the same as for the rotation algebras $A_\theta$. We give a slightly different proof, adapted from Proskurin's (\cite{P00}, \cite{P00b}) classification of the $C^*$-algebras $\Atens$.

\begin{prop}
 For irrational parameters $\theta$, the classification of the $C^*$algebras $\Afrei$ depends on $\theta$ in the same way as it does for the rotation algebras $A_\theta$: 
\[\Afrei \cong \Afreimu \qquad \Longleftrightarrow \qquad \theta = \pm \mu \textnormal{ mod }\Z\]
The same holds true for $\Atens$ and $\Aver$.
\end{prop}
\Proof For $\theta = \pm \mu \textnormal{ mod }\Z$ we have either $\Afrei = \Afreimu$ or $\Afrei \cong \Afreiminus =\Afreimu$ via $u \leftrightarrow v, v\leftrightarrow u$. For the converse, let $\alpha:\Afrei \stackrel{\cong}{\to} \Afreimu$ be an isomorphism between $\Afrei$ and $\Afreimu$, and let $\phi:\Afrei\surj A_\theta$ and $\psi:\Afreimu\surj A_\mu$ be the canonical surjections onto the rotation algebras. 
Then  $\alpha(\ker(\phi)) = \ker(\psi)$ by Proposition \ref{JGroesst}. Thus, $\alpha$ induces an isomorphism from $A_\theta$ to $A_\mu$.

Use Corollary \ref{JGroesstInAtens} for the cases of $\Atens$ and $\Aver$.\qed

\section{A decomposition of the isometries in $\Afrei$} \label{SectDecomp}

\noindent
This section is  the heart of this article, since we will develop the main tools for the investigation of $\Afrei$. We first present a Wold  decomposition of the product $uv$ of the isometries $u$ and $v$ in $\Afrei$ on some Hilbert space $H$. 
This gives rise to a subspace $H_0$ of $H$, on which $uv$ is a unitary. It turns out that $u$ and $v$ are unitaries on $H_0$, too. A study of $u$ and $v$ on the orthogonal complement of $H_0$ reveals essential parts of their structure.

This yields an embedding $\iota$ of $\Afrei$ into $\Asort$ from which we may obtain a lot of information in $\Afrei$. Our approach is adapted from the work of Berger, Coburn and Lebow (\cite{BCL78}), who investigated the $C^*$-algebra generated by two commuting isometries; this is the case $\Afrei$ with $\theta=0$.

Using this embedding $\iota$, we show that the ideal $J=\langle1-uu^*,1-vv^*\rangle\lhd\Afrei$ is isomorphic to $\J$. This will be of crucial use for the computation of the $K$-theory of $\Afrei$. Also, this proves that $\Afrei$ is not exact, because $\Aup$ is not.
Finally, we apply the machinery of this section to the $C^*$-algebra $\Atens$.\\

The well-known Wold decomposition states that every isometry $v$ on a Hilbert space $H$ is of the form $v=v_u\oplus v_s$ for a decomposition  $H=H_u\oplus H_s$ of the Hilbert space. The operator $v_u$ is a unitary and $v_s$ is an (amplified) copy  of the unilateral shift. Hence, the unilateral shift is ``\emph{the}'' isometry and the Wold decomposition reveals the unitary and the shift  part of an isometry.

In their article from 1978 (\cite{BCL78}), Berger, Coburn and Lebow investigated the representation and index theory for $C^*$-algebras generated by commuting isometries. They noticed that the Wold decomposition of the \emph{product} of the commuting isometries yields a subspace for the unitary part, on which  the \emph{single} isometries are unitaries as well. Thereby, they managed to represent the commuting isometries as tensor products of some operators. 
This approach may be slightly modified to the more general case of twisted commuting isometries --  hence to the $C^*$-algebra $\Afrei$. The point is, that the twist of the multiplication does not affect the relevant subspaces.

Let $(e_n)_{n\in\N_0}$ be an orthonormal basis of $\ell^2(\N_0)$,  $S\in\LN$ be the unilateral shift and $d(\lambda)\in\LN$ be the rotation operator, given by $d(\lambda)e_n=\lambda^ne_n$ for $n\in\N_0$ and $\lambda\in\C$ of absolute value one.

\begin{prop}  \label{AfreiDecomp}
 Let $\pi:\Afrei\to\LH$ be a unital representation of $\Afrei$ and let 
 \begin{itemize}
   \item $p:=1-\pi(uvv^*u^*)=1-\pi(vuu^*v^*)\in\LH$ be the defect projection of the isometry $\pi(uv)$,
   \item $H_0$ be the set of all vectors $\xi\in H$ such that for all $n>0$ there exists a $\xi_n\in H$ such that $\xi=\pi(uv)^n\xi_n$,
   \item and let $K$ be the closed linear span of all $\pi(uv)^n\eta$, where $n\geq 0$ and $\eta\in pH$.
 \end{itemize}
Then the following holds:
\begin{itemize}
 \item[(i)] $H_0$ and $K$ are closed linear subspaces of $H$ with $H_0^\perp=K$, and there is a Hilbert space isomorphism $K\cong \ell^2(\N_0)\otimes pH$ via $\pi(uv)^n\eta \leftrightarrow e_n\otimes \eta$.
 \item[(ii)] The restriction $\pi(uv)_{|H_0}\in\LHH{H_0}$ is a unitary, and $\pi(uv)_{|K} \cong S \otimes 1$ under the Hilbert space isomorphism of (i).
 \item[(iii)] The restrictions $\pi(u)_{|H_0}$ and $\pi(v)_{|H_0}$ are unitaries.
 \item[(iv)] Put $\bar u:=u(1-vv^*)+(1-uu^*)v^*$ and $\bar p:=v(1-uu^*)v^*\in\Afrei$. 
Then \linebreak
$\bar u':=\pi(\bar u)$ is a unitary on $pH$ and $\bar p':= \pi(vv^*)p=p\pi(vv^*)=\pi(\bar p)\in\LHH{pH}$ is a projection.
 \item[(v)] The restrictions $\pi(u)_{|K}$ and $\pi(v)_{|K}$ are operators on $K$ and they are of the following form, using the Hilbert space isomorphism of (i):
 \begin{align*}
\pi(u)_{|K}&\cong Sd(\lambda)\otimes \bar u'\bar p' + d(\lambda)\otimes \bar u'(1-\bar p')\\
\pi(v)_{|K}&\cong d(\lambda)^*S\otimes(1-\bar p')\bar u'^* + d(\lambda)^*\otimes \bar p'\bar u'^*
 \end{align*}
\end{itemize}
\end{prop}
\Proof To simplify the notation, we write $w:=\pi(uv)$. 

(i) To see that $H_0$ is closed, check that $(1-w^n(w^*)^n)\xi=0$ for all $n>0$, whenever $\xi\in H$ is a limit of a sequence $(\xi_k)_{k\in\N}\subset H_0$. 
Use $pH \perp wH$ to show $H_0\subset K^\perp$, 
and show inductively $\xi=w^n((w^*)^n\xi)$ for $n>0$ and $\xi\in K^\perp$ to obtain $H_0\supset K^\perp$ (use $\scp{p(w^*)^n\xi}{p(w^*)^n\xi}=0$ for all $n\geq 0$ to show $w^nww^*(w^*)^n\xi=w^n(w^*)^n\xi=\xi$).
If $(\eta_i)_{i\in I}$ forms an orthonormal basis of $pH$ for some index set $I$, the ensemble $(w^n\eta_i)_{n\in \N_0, i\in I}$ forms an orthonormal basis of $K$ as well as the elements $e_n\otimes\eta_i$ for $\ell^2(\N_0)\otimes pH$. This yields the Hilbert space isomorphism $K\cong \ell^2(\N_0)\otimes pH$.

(ii) We have $wH_0 = H_0$ and hence the isometry $w_{|H_0}$ is surjective.
Using the Hilbert space isomorphism of (i), we compute:
\[w(e_n\otimes\eta) \leftrightarrow w(w^n\eta)=w^{n+1}\eta \leftrightarrow e_{n+1}\otimes\eta\]

(iii) Since $\pi(u)$ and $w$ commute up to a scalar, we get $\pi(u)H_0\subset H_0$. Furthermore, $\pi(u)H_0= H_0$, because we can write $\xi\in H_0$ as $\xi=\pi(u)\pi(v)\xi_1$ for some $\xi_1\in H$. Then $\pi(v)\xi_1$ is in $H_0$ as $\pi(v)\xi_1=w^n(\lambda^{-n}\pi(v)\xi_{n+1})$.

(iv) A simple algebraic calculation shows $\bar u\bar u^*=\bar u^*\bar u=1-uvv^*u^*$.

(v) Since $\pi(u)H_0=H_0$, we get that $\pi(u)H_0^\perp\subset H_0^\perp$; likewise  $\pi(v)H_0^\perp\subset H_0^\perp$.

From $\bar u\bar p=(1-uu^*)v^*$ we conclude $\bar u'\bar p' \eta = \pi(v^*)p\eta=\pi(v^*)\eta$ for any $\eta\in pH$. We use this and $\bar u(1-\bar p) =\bar u-\bar u\bar p=u(1-vv^*)$ for the following computation, where $n\in \N_0$ and $\eta\in pH$. Under the Hilbert space isomorphism of (i), we have:
\begin{align*}
 \pi(u)(e_n\otimes\eta) &\leftrightarrow \pi(u)w^n\eta \\
 &=\lambda^nw^n\pi(u)\eta \\
 &=\lambda^n(w^n\pi(u)\pi(vv^*)\eta +w^n\pi(u)(1-\pi(vv^*))\eta) \\
 &=\lambda^n(w^{n+1}\pi(v^*)\eta +w^n\pi(u(1-vv^*))\eta) \\
 &=\lambda^n(w^{n+1}\bar u'\bar p'\eta +w^n\bar u'(1-\bar p')\eta) \\
 &\leftrightarrow (Sd(\lambda)\otimes \bar u'\bar p' + d(\lambda)\otimes\bar u'(1-\bar p'))(e_n\otimes\eta)
\end{align*}
Note, that $\bar u'$ and $\bar p'$ are operators on $pH$, thus $\bar u'\bar p'\eta, \bar u'(1-\bar p')\eta\in pH$ for $\eta\in pH$. 

Analogously, we use $\bar p\bar u^*=v(1-uu^*)$ and $(1-\bar p')\bar u'^*\eta=\pi(u^*)\eta$ for any $\eta\in pH$. This yields:
\begin{align*}
 \pi(v)(e_n\otimes\eta) &\leftrightarrow \pi(v)w^n\eta \\
 &=\lambda^{-n}(w^n\pi(v)\pi(uu^*)\eta +w^n\pi(v)(1-\pi(uu^*))\eta) \\
 &=\lambda^{-n}(\bar \lambda w^{n+1}\pi(u^*)\eta +w^n\pi(v(1-uu^*))\eta) \\
 &=\lambda^{-n}(\bar \lambda w^{n+1}(1-\bar p')\bar u'^*\eta +w^n\bar p'\bar u'^*\eta) \\
 &\leftrightarrow (\bar\lambda Sd(\lambda)^*\otimes(1-\bar p')\bar u'^*+ d(\lambda)^*\otimes\bar p'\bar u'^*)(e_n\otimes\eta)
\end{align*}
Using $\bar\lambda Sd(\lambda)^*=d(\lambda)^*S$, we get the stated result.
\qqed


From Proposition \ref{AfreiDecomp}, we deduce the existence of an embedding $\iota$ of $\Afrei$ into $\Asort$. The unital free product of $\CS$ with $\C^2$ can be described in several ways.
We denote by $\Z_2$ the quotient of $\Z$ by $2\Z$.

\begin{lem} \label{AupGestalt}
 The following $C^*$-algebras are isomorphic:
\begin{itemize}
 \item[(i)] The unital free product $\Aup$ of $\CS$ and $\C^2$.
 \item[(ii)] The universal $C^*$-algebra generated by a unitary $\bar u$ and a projection $\bar p$.
 \item[(iii)] The (full) group $C^*$-algebra $C^*(\Z * \Z_2)$.
 \item[(iv)] The crossed product $C^*(\F_2)\rtimes\Z_2$ under the action that swaps the generators.
\end{itemize}
\end{lem}
\Proof $\C^2$ is the universal unital $C^*$-algebra generated by a projection $\bar p$ under the identification $\left(\begin{array}{c} 1 \\ 0 \end{array}\right)\leftrightarrow\bar p$. As $\CS$ can be viewed as the universal $C^*$-algebra generated by a unitary $\bar u$, this yields the isomorphism of (i) and (ii).

Furthermore, $C^*(\Z)\cong\CS$  and $C^*(\Z_2)\cong\C^2$. Hence the (full) group $C^*$-algebra of $\Z * \Z_2$ may be written as $C^*(\Z * \Z_2)\cong C^*(\Z) *_\C C^*(\Z_2)\cong \Aup$.

From an article by Murphy in 1996 (\cite[proof of th. 6.2]{M96}), we get the isomorphism of (ii) and (iv). The group $C^*$-algebra $C^*(\F_2)$ may be seen as the universal $C^*$-algebra generated by two unitaries $a$ and $b$ without any further relations. The automorphism that swaps $a$ and $b$ is of order two, hence we can form the crossed product with $\Z_2$. We write $C^*(\F_2)\rtimes\Z_2$ as the universal $C^*$-algebra generated by two unitaries $a$ and $b$ together with a symmetry $z$ (i.e. a unitary $z$ with $z^2=1$ or equivalent $z=z^*$) such that $a=zbz$ and $b=zaz$. Since $b$ may be built out of $a$ and $z$, this is just the universal $C^*$-algebra generated by a unitary $a$ and a symmetry $z$.
The homomorphism from $C^*(\bar u, \bar p)$ to $C^*(\F_2)\rtimes\Z_2$, mapping $\bar u\mapsto a$ and $\bar p\mapsto (z+1)/2$ is inverse to the one in the converse direction, mapping $a\mapsto \bar u$ and $z\mapsto 2\bar p -1$. \qqed

From now on, we will view $\Aup$ as the universal $C^*$-algebra generated by a unitary $\bar u$ and a projection $\bar p$.
For the next theorem, note that the generators of $\Afrei$ are denoted by $u$ and $v$ as well as those of $\Aver$ (recall Remark \ref{DefAver}).

\begin{thm} \label{SortEmbed}
There is an embedding of $\Afrei$ into the $C^*$-algebra \linebreak $\Asort$ of the following form.
\begin{align*}
 \iota:\Afrei&\inj\Asort \\
 u&\mapsto vu\otimes \bar u\bar p + u\otimes \bar u(1-\bar p)\\
 v&\mapsto u^*v\otimes (1-\bar p)\bar u^* + u^*\otimes \bar p \bar u^*
\end{align*}
\end{thm}
\Proof The homomorphism $\iota$ exists by the universal property of $\Afrei$. Furthermore, we have a homomorphism $\iota_0:\Afrei\to A_\theta$, mapping $u\mapsto u$ and $v\mapsto v$. The direct sum $\iota_0\oplus\iota$
is injective. Indeed, let $\pi:\Afrei\inj\LH$ be a unital, faithful representation of $\Afrei$. By Proposition \ref{AfreiDecomp} we can decompose $H=H_0\oplus \left(\ell^2(\N_0)\otimes pH\right)$ such that:
\begin{align*}
&\pi(u)=\pi(u)_{|H_0}\oplus \left[Sd(\lambda)\otimes \bar u'\bar p' + d(\lambda)\otimes \bar u'(1-\bar p')\right] \\
&\pi(v)=\pi(v)_{|H_0}\oplus \left[d(\lambda)^*S\otimes(1-\bar p')\bar u'^* + d(\lambda)^*\otimes \bar p'\bar u'^*\right]
\end{align*}
Thus, we have two homomorphism $\sigma_0:A_\theta\to\LHH{H_0}$, $u\mapsto \pi(u)_{|H_0}$ and $v\mapsto \pi(v)_{|H_0}$, and $\sigma: \Asort\to\LHH{\ell^2(\N_0)\otimes pH}$, $u\otimes 1\mapsto d(\lambda)\otimes 1$, $v\otimes 1\mapsto S\otimes 1$, $1\otimes\bar u\mapsto 1\otimes\bar u'$ and $1\otimes\bar p\mapsto1\otimes\bar p'$. Hence $\pi=(\sigma_0\oplus\sigma)\circ(\iota_0\oplus\iota)$, and $\iota_0\oplus\iota$ is injective.

Consider now the homomorphism $\tau:\Asort\to A_\theta$, given by $u\otimes 1\mapsto u$, $v\otimes 1\mapsto uv$, $1\otimes\bar u\mapsto 1$, $1\otimes\bar p\mapsto 0$. Then $\iota_0=\tau\circ\iota$. Thus, $\iota(x)=0$ implies $\iota_0(x)\oplus\iota(x)=0$, which yields $x=0$. Therefore $\iota$ is injective and the proof is complete.

If $\theta$ is irrational, we can simplify the proof using Proposition \ref{JGroesst}. The restriction $\pi_{|H_0}:\Afrei\to\LHH{H_0}$ is isomorphic to the canonical map $\Afrei\to A_\theta$, since $C^*(\pi_{|H_0}(u),\pi_{|H_0}(v))$ is isomorphic to $A_\theta$ ($A_\theta$ is simple). Thus, the kernel of the restriction $\pi_{|H_0^\perp}$ is contained in the kernel of the map $\pi_{|H_0}$ by Proposition \ref{JGroesst}. Hence the restriction $\pi_{|H_0^\perp}=\sigma\circ\iota$ is injective, and therefore also $\iota$. 
\qqed

\begin{rem} \label{ListeIota}
 Direct computations show that we have the following assignment of elements (for $i\in\N_0$) under the map $\iota:\Afrei\to\Asort$.
 \begin{align*}
  u&\mapsto vu\otimes \bar u\bar p + u\otimes \bar u(1-\bar p)\\
  v&\mapsto u^*v\otimes (1-\bar p)\bar u^* + u^*\otimes \bar p \bar u^*\\
  uv&\mapsto v\otimes 1\\
  1-uu^*&\mapsto (1-vv^*)\otimes\bar u\bar p\bar u^*\\
  1-vv^*&\mapsto (1-vv^*)\otimes(1-\bar p)\\
  (uv)^i\bar u&\mapsto v^iu(1-vv^*)\otimes \bar u\\
  (uv)^i\bar u^*&\mapsto v^iu^*(1-vv^*)\otimes \bar u^*\\
  (uv)^i\bar p&\mapsto v^i(1-vv^*)\otimes \bar p
 \end{align*}
Here, $\bar u:= u(1-vv^*)+(1-uu^*)v^*$ and $\bar p:=v(1-uu^*)v^*$ in $\Afrei$ are defined as in Proposition \ref{AfreiDecomp}(iv).

From this embedding $\iota$, we immediately see that the range projections $uu^*$ and $vv^*$ in $\Afrei$ do \emph{not} commute (cf. Remark \ref{ProjKomm}).
\end{rem}

This remark gives us an idea of an explicit picture of the ideal $J$ in $\Afrei$. Note that the elements $\bar u$ resp. $\bar p$ in $J$ are mapped under $\iota$ to $u(1-vv^*)\otimes \bar u$ resp. $(1-vv^*)\otimes \bar p$. This reveals the structure of $\Aup$ in the second component, if we ``untwist'' the multiplication with the unitary $u$ in the image of $\bar u$. In the first component we have a  copy of the algebra of compact operators, if we consider elements $(uv)^i\bar p(uv)^*j$. This yields $J\cong\J$, which will be verified in the next theorem. 

Recall that we see the compact operators as the universal $C^*$-algebra generated by elements $x_i$, $i\in\N_0$ such that $x_i^*x_j = \delta_{ij}x_0$.

\begin{thm}\label{JGestalt}
 The defect ideal $J=\langle1-uu^*,1-vv^*\rangle\lhd\Afrei$ is isomorphic to $\J$ via $(uv)^i\bar u\leftrightarrow\bar u\otimes x_i$ and $(uv)^i\bar p\leftrightarrow\bar p\otimes x_i$ for $i\in\N_0$.
\end{thm}
\Proof Like in Proposition \ref{DefIdealAverAtens}, we first use a twisted version of $\J$ and show that it is isomorphic to $\iota(J)\cong J$. We then apply Lemma  \ref{TensKompUnTwist} to get the isomorphism to $\J$.

Let $\alpha_i$ defined by $\alpha_i(\bar u):=\lambda^i\bar u$ and $\alpha_i(\bar p):=\bar p$ be automorphisms on $\Aup$ for $i\in\N_0$ and form the twisted tensor product with the compacts, $(\Aup)\otimes_\alpha\K$ as defined in Definition \ref{DefTwistKomp}.
 Hence  $(\Aup)\otimes_\alpha\K$ is an ideal in the universal $C^*$-algebra $B_\alpha$ which in turn is generated by a unitary $\bar u$, a projection $\bar p$, elements $x_i$ for $i\in\N_0$, and the relations $x_i^*x_j=\delta_{ij}x_0$, $\bar ux_i=\lambda^ix_i\bar u$, and $\bar px_i=x_i\bar p$.
By the universal property, we obtain the following homomorphism: 
\begin{align*}
&\phi:B_\alpha\to\Asort\\
 & \bar u\mapsto u\otimes\bar u, \quad \bar p\mapsto 1\otimes \bar p, \quad x_i\mapsto v^i(1-vv^*)\otimes 1
\end{align*}

The ideal $(\Aup)\otimes_\alpha\K$ is spanned by all elements $\rho x_ix_j^*$, where $\rho$ is a $^*$-monomial in $\bar u$ and $\bar p$, and $i,j\in\N_0$.
They are mapped by $\phi$ to $u^kv^i(1-vv^*)(v^*)^j\otimes\rho$ with $\sigma(\rho)= u^k$, for some $k\in\Z$. Here, $\sigma:\Aup\to\CS$ is given by $\bar u\mapsto u$, $\bar p\mapsto 1$, and we use the canonical embedding of $\Aup$ into $B_\alpha$. 
The image of $J=\langle 1-uu^*,1-vv^*\rangle\lhd\Afrei$ under $\iota$ is exactly the closed linear span of these elements $u^kv^i(1-vv^*)(v^*)^j\otimes\rho\in\Asort$. Indeed, let $\rho'\in J$ be the $^*$-monomial in $\bar u$ and $\bar p\in\Afrei$ analogous to $\rho\in\Aup$. Then $\iota((uv)^i\rho'((uv)^*)^j)$ equals $u^kv^i(1-vv^*)(v^*)^j\otimes\rho$ up to a scalar, because $\iota(\bar u)=u(1-vv^*)\otimes \bar u$,  $\iota(\bar p)=(1-vv^*)\otimes \bar p$ and $\iota(uv)=v\otimes1$  (cf. Remark \ref{ListeIota}). On the other hand, the closed linear span by the elements $u^kv^i(1-vv^*)(v^*)^j\otimes\rho$ is an ideal in $\iota(\Afrei)$ containing $\iota(1-uu^*)$ and $\iota(1-vv^*)$.
We conclude that the image of $(\Aup)\otimes_\alpha\K$ under $\phi$ is exactly $\iota(J)$.

It remains to show, that the restriction of $\phi$ to $(\Aup)\otimes_\alpha\K$ is injective. 
%
%
For this, let $\pi:B_\alpha\inj\LH$ be a unital, faithful representation of $B_\alpha$ and decompose the Hilbert space $H=\left(\bigoplus_{i\in\N_0}\pi(x_ix_i^*)H\right)\oplus K$, where $K\subset H$ is a subspace of $H$. There is a Hilbert space isomorphism $\bigoplus_{i\in\N_0}\pi(x_ix_i^*)H\cong\ell^2(\N_0)\otimes H_0$, where $H_0:=\pi(x_0)H$, via $\pi(x_i)\eta\leftrightarrow e_i\otimes\eta$ (see the proof of Proposition \ref{DefIdealAverAtens}).
 We  now compute the form of the elements of $B_\alpha$ on the space $K^\perp$. Note, that $\pi(\bar u)H_0=H_0$, since $\bar ux_0=x_0\bar u$ and $\bar u$ is a unitary. Hence $\tilde u:=\pi(\bar u)_{|H_0}\in\LHH{H_0}$ is a unitary on $H_0$. Also $\tilde p:=\pi(\bar p)_{|H_0}\in\LHH{H_0}$ acts on $H_0$, because $\bar px_0=x_0\bar p$.
On $K^\perp$, the operators have the following from, using the Hilbert space isomorphism $K^\perp\cong\ell^2(\N_0)\otimes H_0$.
\begin{align*}
 &\pi(\bar u)_{|K^\perp}(e_i\otimes\eta)  &&\leftrightarrow \pi(\bar ux_i)\eta=\lambda^i\pi(x_i)\pi(\bar u)\eta   &&\leftrightarrow (d(\lambda)\otimes\tilde u)(e_i\otimes\eta)\\
&\pi(\bar p)_{|K^\perp}(e_i\otimes\eta) &&\leftrightarrow \pi(\bar px_i)\eta =\pi(x_i)\pi(\bar p)\eta   &&\leftrightarrow (1\otimes\tilde p)(e_i\otimes\eta)\\
&\pi(x_j)_{|K^\perp}(e_i\otimes\eta)  &&\leftrightarrow \pi(x_jx_i)\eta  =\delta_{i0}\pi(x_j)\eta   &&\leftrightarrow (S^j(1-SS^*)\otimes 1)(e_i\otimes\eta)
\end{align*}
Here, we used $x_jx_i=x_jx_0^*x_i=\delta_{i0}x_j$.
Because of $\pi(x_j)_{|K}=0$ for all $j\in\N_0$ and as $(\Aup)\otimes_\alpha\K$ is spanned by elements $\rho x_ix_j^*$, the restriction of $\pi$ to $(\Aup)\otimes_\alpha\K$ is a unital, faithful map to $\LHH{K^\perp}$ mapping
$\bar u\mapsto d(\lambda)\otimes\tilde u$, $\bar p\mapsto 1\otimes \tilde p$, and $x_j\mapsto (S^j(1-SS^*))\otimes 1$.

On the space $K^\perp\cong\ell^2(\N_0)\otimes H_0$, the operators $d(\lambda)\otimes1$, $S\otimes1$, $1\otimes \tilde u$ and $1\otimes \tilde p$ give rise to a representation of $\Asort$. Thus we get a homomorphism: 
\begin{align*}
 &\tau:\Asort\to\LHH{K^\perp}\\
 &u\otimes 1\mapsto d(\lambda)\otimes 1, \quad v\otimes 1\mapsto S\otimes 1,\\ 
 &1\otimes\bar u\mapsto 1\otimes\tilde u,\quad 1\otimes\bar p\mapsto 1\otimes \tilde p
\end{align*}
We conclude that $\pi_{|K^\perp}=\tau\circ\phi$ on $(\Aup)\otimes_\alpha\K$, which proves injectivity for $\phi$ on $(\Aup)\otimes_\alpha\K$.
Therefore, $(\Aup)\otimes_\alpha\K$ is isomorphic to $\iota(J)\cong J$ and we may apply Lemma \ref{TensKompUnTwist} to complete the proof.
\qqed

In the case $\theta=0$, Murphy (\cite[proof of th. 6.2]{M96}) obtained the same picture for $J\lhd\Afrei$ by similar means. But the case of general $\theta$ remained untouched.

This picture of $J$ enables us to compute its $K$-theory, since the $K$-groups of free products of $C^*$-algebras are completely understood by the work of Cuntz (\cite{C82}). We then may derive the $K$-groups of $\Afrei$ for all $\theta$. This is done in section \ref{SectKTheorie}.

Furthermore, Theorem \ref{JGestalt} reveals that $\Afrei$ is not exact, since $\Aup$ is not. The latter may be deduced from Lemma \ref{AupGestalt}.




\begin{prop}\label{AupNotExact}
 The $C^*$-algebra $\Aup$ is not exact.
\end{prop}
\Proof 
By Lemma \ref{AupGestalt}, $\Aup$ is isomorphic to $C^*(\F_2)\rtimes\Z_2$. Hence it contains $C^*(\F_2)$ as a $C^*$-subalgebra. 
Since $C^*(\F_2)$ is not exact due to Simon Wassermann (see \cite{W90} or \cite{W78} and \cite{W76}) and since a $C^*$-subalgebra of an exact $C^*$-algebra is exact again (a $C^*$-algebra is exact if and only if it may be embedded into $\On{2}$, \cite{KP00}), we conclude that $\Aup$ cannot be exact. \qqed

\begin{thm}
 $\Afrei$ is not exact.
\end{thm}
\Proof We have $\Aup\subset\J\cong J\subset\Afrei$. \qqed

\begin{rem} \label{WeiteresAfrei} 
The following is an overview on some further properties of $\Afrei$.
\begin{itemize}
 \item[(a)] As $C^*(\F_2)$ embeds into $\Afrei$ and as the first $C^*$-algebra is \emph{not locally reflexive}, nor is the second. (cf. \cite[Cor. 9.1.6 and Lem. 9.2.8]{BO08}).
 \item[(b)] Murphy (\cite[Th. 3.3]{M03}) showed that the (full) $C^*$-algebra of the free product of a non-trivial, countable, free group $F$ and a non-trivial, countable, amenable group $Z$ is primitive. Hence $\Aup$ is primitive by Lemma \ref{AupGestalt}. 
 Now let $I_1, I_2\lhdnicht\Afrei$ be two non-zero ideals in $\Afrei$. For irrational parameter $\theta$, we have $I_1, I_2\subset J$ by Proposition \ref{JGroesst}. Using the isomorphism of Theorem \ref{JGestalt}, we see that these ideals are isomorphic to $I_1'\otimes\K$ and $I_2'\otimes\K$ for some non-zero ideals  $I_1'$ and $I_2'$ in $\Aup$. 
 As $\Aup$ is not primitive, $I_1'$ and $I_2'$ have non-zero intersection, thus the same holds for $I_1$ and $I_2$. We infer that $\Afrei$ is \emph{primitive} for irrational $\theta$.
 \item[(c)] The $C^*$-algebra $\Aup$ is not only a $C^*$-subalgebra of $\Afrei$ but also a \emph{hereditary} $C^*$-subalgebra. Under the isomorphism of Theorem \ref{JGestalt}, we have that $\Aup$ is isomorphic to the $C^*$-subalgebra of $\Afrei$ generated by the elements $\bar u:= u(1-vv^*)+(1-uu^*)v^*$ and $\bar p:=v(1-uu^*)v^*$ (cf. Remark \ref{ListeIota} and Theorem \ref{JGestalt} and using the embedding of $\Aup$ into $\Aup\otimes\K$ by $x\mapsto x\otimes x_0$). This in turn coincides with the compression $p\Afrei p$ of $\Afrei$ by the projection $p:= 1-uvv^*u^*$ (which is the unit of $C^*(\bar u, \bar p)$). Indeed, the projection $p$ corresponds to the minimal projection $1\otimes x_0$ in $\J$, thus:
 \[\Aup\cong(1\otimes x_0)(\J)(1\otimes x_0)\cong  p J  p =p\Afrei p\]
\end{itemize}
\end{rem}

\begin{rem}\label{EinbettungFallAtens}
The techniques developed in this section may also be applied to  $\Atens$, although this  $C^*$-algebra is much less complicated. Nevertheless, let us quickly review this case.
\begin{itemize}
\item[(a)] In analogy to Proposition \ref{AfreiDecomp}, we consider the Wold decomposition of the isometry $\pi(uv)$, if $\pi$ is a unital representation of $\Atens$. In this case, the unitaries $\bar u'$ and $\bar p'$ in $\LHH{pH}$ fulfill the relations $\bar p'\bar u'=\bar p'\bar u'\bar p'$ (with the same notations as in Prop. \ref{AfreiDecomp}). Accordingly, we obtain an embedding 
\[\iota':\Atens\inj\AsortT. \]
Here, $L$ is the quotient of $\Aup$ by the ideal generated by 
$\bar p\bar u-\bar p\bar u\bar p$. We deduce that  the defect ideal $\langle1-uu^*,1-vv^*\rangle\lhd\Atens$ is isomorphic to $L\otimes\K$.
\item[(b)] The $C^*$-algebra $L$ is much smaller than $\Aup$. In fact, there is an embedding $\psi$ of $L$ into $ M_2(\Toepl)=M_2(\C)\otimes \Toepl$, by  
\[
\bar u\mapsto 
\begin{pmatrix}
v^*      & 0 \\
1-vv^*   & v \\
\end{pmatrix}
\qquad \textnormal{and} \qquad
\bar p\mapsto 
\begin{pmatrix}
1      & 0 \\
0   & 0 \\
\end{pmatrix}.
\]
\item[(c)] The following sequence is exact, hence $L$ is nuclear.
\[0\to\K\to L \to \CS\oplus\CS\to 0\]
Tensoring this sequence with the compacts $\K$ (and using item (a)) yields a sequence which Kabluchko proved in 2001 (\cite{K01}) to be exact:
\[0\to\K\to J' \to \K\otimes(\CS\oplus\CS)\to 0\]
Here $J'=\langle1-uu^*,1-vv^*\rangle\lhd\Atens$. For this, he showed that the ideals $\langle 1-uu^*\rangle$ and  $\langle 1-vv^*\rangle\lhd\Atens$ are isomorphic to $\Toepl\otimes\K$ (cf. Proposition \ref{DefIdealAverAtens} of our article) and $\langle 1-uu^*\rangle\cap\langle1-vv^*\rangle=\K$ whereas $\langle1-uu^*\rangle\cup\langle1-vv^*\rangle=J'$. Hence $J'/\K\cong (\langle1-uu^*\rangle/\K)\oplus(\langle1-vv^*\rangle/\K)=(\CS\otimes\K)\oplus(\CS\otimes\K)$.
\end{itemize}
\end{rem}

\Proof (b) The projections  $q_0:=\bar p-\bar u^*\bar p\bar u$ and $q_n:=\bar u^nq_0\bar u^{-n}$ for $n\in\Z$ are mutually orthogonal. Furthermore, $\bar uq_n=q_{n+1}\bar u$ for all $n\in\Z$. Thus $\pi(\bar u)$ acts as a bilateral shift on 
$K_0:=\bigoplus_{n\in\Z}\pi(q_n)H\subset H$, if $\pi:L\to\LH$ is a representation of $L$, whereas $\pi_{|K_0}(\bar p)$ projects onto $q_-H:=\sum_{n\leq 0}\pi(q_n)H$.
This yields a representation of $M_2(\Toepl)$ on $\mathcal L(K_0)$ which proves that $\psi$ is injective, if $\pi$ is faithful.

(c) Consider the homomorphism $\sigma:L\to\CS\oplus\CS$, mapping $\bar u\mapsto(u,u)$ and $\bar p\mapsto(1,0)$. Its kernel is the ideal $\langle \bar u\bar p-\bar p\bar u\rangle=\langle q_0\rangle$, where $q_0$ is the projection given by $q_0:=\bar p-\bar u^*\bar p\bar u\in L$. 
This ideal is spanned by matrix units $\bar e_{kl}:=\bar u^k q_0 (\bar u^*)^l$ for $k,l\in\Z$, hence it is isomorphic to the compacts $\K$.
\qed

\section{The ideal structure of $\Afrei$} \label{SectIdStr}

\noindent
For irrational parameters $\theta$, a combination of the results of Proposition \ref{JGroesst} and Theorem \ref{JGestalt} gives a good description of the ideal structure of $\Afrei$. Since every nontrivial ideal $I$ in $\Afrei$ is contained in the ideal $J=\Iuv$ which in turn is isomorphic to $\J$, we conclude that we have a one-to-one correspondence of ideals in $\Afrei$ with ideals in $\Aup$, if $\theta$ is irrational. (Note that every ideal in $A\otimes\K$ is of the form $I'\otimes\K$, if $A$ is unital.) We investigate this correspondence in an explicit way. Furthermore, we prove that not only $\Afrei$ is not exact, but neither is its quotient by the ideal generated by the commutator $[uu^*,vv^*]$.\\

Recall that the $C^*$-subalgebra generated by $\bar u:= u(1-vv^*)+(1-uu^*)v^*$ and $\bar p:=v(1-uu^*)v^*$ in $\Afrei$ is isomorphic to $\Aup$ (see Remark \ref{WeiteresAfrei}(c)).

\begin{lem} \label{DualityTheorem}
Let $\theta$ be irrational.
 Let $I\lhdnicht\Afrei$ and $I'\lhd\Aup$ be ideals, and let $\phi$ resp. $\beta$ be the according quotient maps. Let $\phi(I')=0$ under the identification $I'\subset\Aup\cong C^*(\bar u, \bar p)\subset \Afrei$ (cf. Remark \ref{WeiteresAfrei}(c)).
 Furthermore, let $\psi:\Afrei / I\to (\Aver)\otimes (\Aup / I')$ be a map such that the following diagram commutes (where $\iota$ is the embedding of Theorem \ref{SortEmbed}).
\begin{align*}
 &&\Afrei &&\stackrel{\iota}{\inj} && \Asort \\
 &&\phi \downarrow &&  && \downarrow \textnormal{id}\otimes \beta \\
 &&\Afrei/I &&\stackrel{\psi}{\longrightarrow} && (\Aver) \otimes (\Aup/I')
\end{align*}
 Then $I\cong I'\otimes\K$. In other words, $I$ corresponds to $I'$ under the isomorphism of $J$ and $\J$.
\end{lem}

\Proof By Proposition \ref{JGroesst}, $I$ is an ideal in $J$ which in turn is isomorphic to $\J$. Thus, $I$ is isomorphic to $I_0\otimes \K$ where $I_0$ is an ideal in $\Aup$. Furthermore, we can identify $I\cap C^*(\bar u, \bar p)$ in $\Afrei$ with $I_0$ in $\Aup$. Indeed, under the identification of $C^*(\bar u, \bar p)\subset \Afrei$ and $\Aup$ via $\bar u\leftrightarrow \bar u$ and $\bar p\leftrightarrow \bar p$ (by Remark \ref{WeiteresAfrei}(c)), we see that the isomorphism of $J$ and $\J$ is of the following form (cf. Theorem \ref{JGestalt}):
\[(uv)^ix \longleftrightarrow x\otimes x_i \qquad \forall x\in\Aup\cong C^*(\bar u, \bar p)\subset \Afrei, i\in\N_0\]
Now, an element  $x\in I\cap C^*(\bar u, \bar p)$ is mapped to $x\otimes x_0\in (I_0\otimes \K)\cap ((\Aup)\otimes x_0)$ under the isomorphism of $J$ and $\J$. Thus it is mapped to $x$ in $I_0$, since $I_0$ embeds into $I_0\otimes\K$ by $a\mapsto a\otimes x_0$ (recall that $x_0$ is a minimal projection).

To prove $I'\subset I_0$, recall that $\phi(x)=0$ for all $x\in I'$ by assumption. Hence $x$ is in $I\cap C^*(\bar u, \bar p)=I_0$.
For the converse direction, consider the representation $\pi:\Aver\to\LN$ defined by $\pi(u)=d(\lambda)$ and $\pi(v)=S$. Here $d(\lambda)$ is given by $d(\lambda)e_n=\lambda^ne_n$ for all $n\in\N_0$ and $S$ is the unilateral shift.
Consider now  $\gamma:\Aup\to\LN\otimes(\Aup/I')$, given by $\gamma:=(\pi\otimes\textnormal{id})\circ\psi\circ\phi$. (Again we identify $\Aup=C^*(\bar u, \bar p)\subset\Afrei$.)
Then $\bar u$ is mapped by $\psi\circ\phi=(\id\otimes\beta)\circ\iota$  to $u(1-vv^*)\otimes\beta(\bar u)$ (cf. Remark \ref{ListeIota}). This is mapped to $d(\lambda)(1-SS^*)\otimes\beta(\bar u)$ by  $(\pi\otimes\textnormal{id})$. Since $d(\lambda)(1-SS^*)$ coincides with $1-SS^*$ we infer $\gamma(\bar u)=(1-SS^*)\otimes\beta(\bar u)$.
Similarly, we see $\gamma(\bar p)=(1-SS^*)\otimes\beta(\bar p)$, and we conclude
$\gamma(x)=(1-SS^*)\otimes\beta(x)$ for all $x\in\Aup$ since  $(1-SS^*)\otimes\beta$ is a homomorphism.
Finally, let $x$ be in $I_0=I\cap C^*(\bar u, \bar p)$. Then $\phi(x)=0$ which implies $(1-SS^*)\otimes\beta(x)=\gamma(x)=0$, hence $\beta(x)=0$. Thus $I_0\subset I'$ and we conclude $I_0=I'$.\qqed

We now give a description of some corresponding ideals in the most interesting cases.
A quick look at the $C^*$-algebra $\Aup$ -- again viewed as the universal $C^*$-algebra generated by a unitary $\bar u$ and a projection $\bar p$ -- shows that we have the following canonical maps and the according kernels. First, we have the quotient map to the $C^*$-algebra $L$ (as defined in Remark \ref{EinbettungFallAtens}). Secondly, we can consider the quotient of $\Aup$ such that $\bar u$ and $\bar p$ commute. This is isomorphic to the direct sum $\CS\oplus\CS$ under the assignment $\bar u\leftrightarrow (u,u)$ and $\bar p\leftrightarrow (1,0)$.
 Thirdly, the quotients where $\bar p=1$ or $\bar p=0$ seem to be natural. Hence we get the following maps and ideals:
\begin{itemize}
 \item $\Aup\to L$, mapping $\bar u\mapsto \bar u$ and $\bar p\mapsto \bar p$. The kernel is the ideal generated by $\bar p\bar u-\bar p\bar u\bar p$. Similarly, we consider the ideal generated by $\bar u\bar p-\bar p\bar u\bar p$.
 \item $\Aup\to\CS\oplus\CS$, mapping $\bar u\mapsto (u,u)$ and $\bar p\mapsto (1,0)$. The according ideal is generated by $\bar u\bar p-\bar p\bar u$.
 \item $\Aup \to\CS$, mapping $\bar u\mapsto u$ and $\bar p\mapsto 1$ resp.  $\bar u\mapsto u$ and $\bar p\mapsto 0$. The ideals are generated by $1-\bar p$ resp. $\bar p$.
\end{itemize}
These ideals are related in an obvious way by: 
\begin{itemize}
 \item $\langle \bar p\rangle + \langle1-\bar p\rangle = \Aup$
 \item $\langle \bar p\rangle \cap \langle1-\bar p\rangle=\langle \bar u\bar p-\bar p\bar u\rangle$ (Write the above map $\Aup\to\CS\oplus\CS$ as the direct sum of the two above maps $\Aup \to\CS$.)
 \item $\langle \bar p\bar u(1-\bar p)\rangle + \langle (1-\bar p)\bar u\bar p\rangle = \langle \bar u\bar p-\bar p\bar u\rangle$
\end{itemize}
In the next two lemmas we work out the connection between the relations on $u$ and $v$ in $\Afrei$ on the one side, and $\bar u$ and $\bar p$ in $\Aup$ on the other side.

\begin{lem} \label{LemKorresp}
 Let $A$ be a unital $C^*$-algebra and let $u$ and $v$ be isometries in $A$ such that $uv=\lambda vu$ where $\lambda\in\C$ is of absolute value one. Put $\bar u:=u(1-vv^*)+(1-uu^*)v^*$, $\bar p:= v(1-uu^*)v^*$ and $p:=1-vuu^*v^*$ as in Proposition \ref{AfreiDecomp}(iv).
Then the following equivalences hold.
\begin{itemize}
 \item[(a)] $u$ is unitary $\eq$ $\bar p=0$
 \item[]    $v$ is unitary $\eq$ $p(1 -\bar p)=0$ 
 \item[]    (Note that $p$ is the unit of $C^*(\bar u,\bar p)\subset A$, cf. Remark \ref{WeiteresAfrei}(c).)
 \item[(b)] $u^*v=\bar \lambda vu^* \eq (1-uu^*)v^*u(1-vv^*)=0 \eq \bar p\bar u(1-\bar p)=0$
 \item[]    $(1-uu^*)(1-vv^*)=0 \eq (1-\bar p)\bar u\bar p=0$
 \item[(c)] 
 $u^*v=\bar\lambda vu^*$ and $(1-uu^*)(1-vv^*)=0 \eq
 \bar u\bar p=\bar p\bar u$
\end{itemize}
\end{lem}
\Proof (a)  The second equivalence is due to the equality $p(1-\bar p)=1-vv^*$.

(b) If $(1-uu^*)v^*u(1-vv^*)=0$, then 
\[v^*u=((uu^*)+(1-uu^*))v^*u((vv^*)+(1-vv^*))=\lambda uv^*\]
Secondly $\bar p\bar u(1-\bar p)=v(1-uu^*)v^*u(1-vv^*)$. 

Thirdly $(1-\bar p)\bar u\bar p=(1-vv^*)(1-uu^*)v^*$.\qqed

If $I$ is an ideal in $\Aup$, we have a map from $\Afrei$ to $(\Aver)\otimes(\Aup/I)$, analogous to the embedding  $\iota$ of $\Afrei$ into $\Asort$ (cf. Theorem \ref{SortEmbed}). In the next lemma, we analyze the structure of the elements $u$ and $v$ under this map.

\begin{lem} \label{LemKorresp2}
 Let $I\lhd\Aup$ be an ideal and let $\beta$ denote its quotient map. Define elements in  $(\Aver)\otimes (\Aup / I)$ by
\begin{align*}
 &\tilde u:= vu\otimes \beta(\bar u\bar p) + u\otimes \beta(\bar u(1-\bar p))\\
 &\tilde v:= u^*v\otimes \beta((1-\bar p)\bar u^*) + u^*\otimes \beta(\bar p \bar u^*)
\end{align*}
Then $\tilde u$ and $\tilde v$ are isometries and $\tilde u\tilde v=\lambda\tilde v\tilde u$. Furthermore:
\begin{itemize}
 \item[(a)] If $I=\langle\bar p\rangle$, then $\tilde u$ is even a unitary.
\item[]    If $I=\langle1-\bar p\rangle$, then $\tilde v$ is even a unitary.
 \item[(b)] If $I=\langle\bar p\bar u(1-\bar p)\rangle$, then $\tilde u^*\tilde v=\bar\lambda\tilde v\tilde u^*$.
\item[]    If $I=\langle(1-\bar p)\bar u\bar p\rangle$, then $(1-\tilde u\tilde u^*)(1-\tilde v\tilde v^*)=0$.
\end{itemize}
\end{lem}
\Proof If $I=0$, then $\tilde u$ and $\tilde v$ are isometries with $\tilde u\tilde v=\lambda\tilde v\tilde u$, hence  this also holds in any quotient $(\Aver)\otimes (\Aup / I)$.
The remaining statements can be checked by a straight forward computation. (Cf. also Remark \ref{ListeIota}.)
\qqed

\begin{prop}\label{propKonkreteKorrespondenzen}
 For irrational parameter $\theta$, there are the following correspondences between ideals  in $\Afrei$ and in $\Aup$ under the isomorphism of $J$ and $\J$ of Theorem 
\ref{JGestalt}.
\begin{itemize}
 \item[(a)] $\langle1-uu^*\rangle\lhd\Afrei$ corresponds to $\langle\bar p\rangle\lhd\Aup$.
 \item[]    $\langle1-vv^*\rangle\lhd\Afrei$ corresponds to $\langle1-\bar p\rangle\lhd\Aup$.
\item[(b)] $\langle u^*v-\bar\lambda vu^*\rangle\lhd\Afrei$ corresponds to $\langle\bar p\bar u(1-\bar p)\rangle\lhd\Aup$.
 \item[]   $\langle (1-uu^*)(1-vv^*)\rangle\lhd\Afrei$ corresponds to $\langle(1-\bar p)\bar u\bar p\rangle\lhd\Aup$.
\item[(c)] $\langle u^*v-\bar\lambda vu^*, (1-uu^*)(1-vv^*)\rangle\lhd\Afrei$ corresponds to $\langle \bar u\bar p-\bar p\bar u\rangle\lhd\Aup$.
(Recall that a description of the quotient of $\Afrei$ by this ideal is given in Proposition \ref{AtensUndAver}.)
\end{itemize}
\end{prop}

\Proof (a) Under the quotient map $\phi:\Afrei\to\Afrei / \langle1-uu^*\rangle$, we have $\phi(\bar p)=0$ by Lemma \ref{LemKorresp} since $\phi(u)$ is a unitary. 
Thus, $\phi(\langle\bar p\rangle)=0$ in the sense of Lemma \ref{DualityTheorem} (under the identification of $\Aup$ and $C^*(\bar u,\bar p)\subset\Afrei$).
 By Lemma \ref{LemKorresp2}, there is a map  $\psi:\Afrei / \langle1-uu^*\rangle\to (\Aver)\otimes (\Aup / \langle\bar p\rangle)$, mapping $\phi(u)\mapsto u\otimes\beta(\bar u)$ and $\phi(v)\mapsto u^*v\otimes \beta(\bar u^*)$, such that the diagram of of Lemma \ref{DualityTheorem} commutes. Thus, $\langle 1-uu^*\rangle=\ker\phi\lhd\Afrei$ corresponds to $\langle\bar p\rangle\lhd\Aup$.

The other correspondences are obtained in exactly the same way.
\qqed

We conclude that the correspondence of ideals in $\Afrei$ and in $\Aup$ matches natural ideal to natural ideals.

In Remark \ref{ProjKomm}, the commutator $[uu^*,vv^*]$ in $\Afrei$ has been considered. We are now able to show that the  quotient of $\Afrei$ by the ideal generated by it is not exact. In some sense, this shows that $\Atens$ and $\Afrei$ differ much more than just by the structure of their range projections. This answers a question of section \ref{SectTwist}.

\begin{lem}\label{LemKorresp4}
  Let $A$ be a unital $C^*$-algebra and let $u$ and $v$ be isometries in $A$ such that $uv=\lambda vu$ for a given $\lambda\in\C$ of absolute value one. With $\bar u\in A$ and $\bar p\in A$ as in Lemma \ref{LemKorresp}, the following relations are equivalent.
\[ [\bar p,\bar u\bar p\bar u^*]=0 \eq [uu^*,vv^*]=0\]
In fact, even the equality $[\bar p,\bar u\bar p\bar u^*]=[uu^*,vv^*]$ holds.
\end{lem}
\Proof A simple algebraic computation shows $\bar p\bar u\bar p\bar u^*=vv^*(1-uu^*)$.
\qqed

\begin{lem} \label{LemKorresp5}
Consider $I=\langle [\bar p,\bar u\bar p\bar u^*]\rangle\lhd\Aup$ and let $\tilde u$ and $\tilde v$ be the isometries in $(\Aver)\otimes (\Aup / I')$ defined as in Lemma \ref{LemKorresp2}. 
Then $[\tilde u\tilde u^*,\tilde v\tilde v^*]=0$.
\end{lem}
\Proof Cf. Remark \ref{ListeIota}. \qqed

\begin{prop}\label{KorAfreikom}
For irrational $\theta$, the ideal $\langle [uu^*,vv^*]\rangle\lhd\Afrei$ corresponds to the ideal $\langle [\bar p, \bar u\bar p\bar u^*]\rangle\lhd\Aup$ under the isomorphism of $J$ and $\J$. 
\end{prop}
\Proof The proof is analogous to that of Proposition \ref{propKonkreteKorrespondenzen}.\qqed

Next, we show that a quotient of $\Aup$ by $\langle [\bar p,\bar u\bar p\bar u^*]\rangle$ contains $\Aup$ as a $C^*$-subalgebra, hence it cannot be exact. This transfers directly to the quotient of $\Afrei$ by the commutator $[uu^*,vv^*]$.

\begin{prop}\label{Afreikomnichtexakt} \label{Kommutator}
\begin{itemize}
 \item[(a)] The quotient $D$ of $\Aup$ by $\langle [\bar p,\bar u\bar p\bar u^*]\rangle$ is not exact.
 \item[(b)] The quotient $B$ of $\Afrei$ by $\langle [uu^*,vv^*]\rangle$ is not exact, if $\theta$ is irrational.
\end{itemize}
\end{prop}
\Proof (a) Consider the elements  
$\bar u':=\begin{pmatrix}
0      & \bar u\\
1   & 0 \\
\end{pmatrix}$ 
and 
$\bar p':=\begin{pmatrix}
\bar p      & 0\\
0   & 0 \\
\end{pmatrix}$ 
in the $2\times 2$-\linebreak matrices $M_2(\Aup)$ over $\Aup$. Then $\bar u'$ is a unitary and $\bar p'$ is a projection. Furthermore, 
$\bar u'\bar p'\bar u'^*=\begin{pmatrix}
0      & 0\\
0   & \bar p \\
\end{pmatrix}$, thus $[\bar p',\bar u'\bar p'\bar u'^*]=0$, which yields a surjective homomorphism from $D$ to $C^*(\bar u', \bar p')\subset M_2(\Aup)$.

But the $C^*$-subalgebra $C^*(\bar u'^2, \bar p')$ of $C^*(\bar u', \bar p')$ is isomorphic to $\Aup$. Indeed, as 
$\bar u'^2=\begin{pmatrix}
\bar u      & 0\\
0   & \bar u \\
\end{pmatrix}$, the homomorphism $\Aup\to C^*(\bar u'^2, \bar p')\subset M_2(\Aup)$ given by $\bar u\mapsto \bar u'^2$, $\bar p\mapsto \bar p'$ may be seen as the map $\id\oplus\rho$. Here,  the homomorphism $\rho:\Aup \to\Aup$ is given by  $\bar u\mapsto \bar u$ and $\bar p\mapsto 0$. 

Thus $C^*(\bar u', \bar p')\subset M_2(\Aup)$ contains the non-exact $C^*$-subalgebra $C^*(\bar u'^2, \bar p')$ and therefore is is not exact.
Since quotients of exact $C^*$-algebras are exact (cf. \cite[Cor. 9.4.3]{BO08}), $D$ cannot be exact.

(b) Let $I$ be a nontrivial ideal in $\Afrei$ and let $\phi$ be the correponding quotient map. Let $I'$ be the ideal in $\Aup$ corresponding to $I$. 
 Consider the following diagram of exact sequences:
\begin{align*}
0 &&\to &&I &&\to &&\Afrei &&\to &&\Afrei/I &&\to &&0 \\
     &&  &&=  &&  &&\bigtriangledown && &&\bigtriangledown \\
0 &&\to &&I &&\to &&J &&\to &&\phi(J) &&\to &&0 
\end{align*}
Compare the lower one with the following exact sequence:
\[0\to I' \otimes \K\to\J \to (\Aup/I')\otimes\K\to 0\]
As $I\cong I' \otimes \K$  and $J\cong \J$, we may conclude that $\phi(J)$ is isomorphic to $(\Aup/I')\otimes\K$.

With $I=\langle [uu^*,vv^*]\rangle$ and $I'=\langle [\bar p,\bar u\bar p\bar u^*]\rangle$, we see that $D$ is a non-exact $C^*$-subalgebra of $B$.\qqed

\section{The $K$-groups of $\Afrei$ and $\Atens$} \label{SectKTheorie}

\noindent
 For the classification of $C^*$-algebras the $K$-theory is the most important ingredient. This section is devoted to the computation of the $K$-groups of $\Atens$ and $\Afrei$. While the case of $\Atens$ is not too difficult to treat, the case of $\Afrei$ relies essentially on the isomorphism of the ideal $J=\Iuv$ with $\J$ found in section \ref{SectDecomp}. By this we can compute its $K$-theory using a general result by Cuntz (\cite{C82}) on the $K$-theory of free products of $C^*$-algebras. 

Consider the following two exact sequences:
\begin{align*}
 0&&\to &&I &&\to &&\Afrei &&\to &&\Aver &&\to &&0 \\
    &&  &&=  &&  &&\bigtriangledown && &&\bigtriangledown \\
 0&&\to &&I &&\to &&J &&\to &&\langle 1-vv^*\rangle &&\to &&0 
\end{align*}
By the lower one, we are able to compute the $K$-groups of the ideal $I=\Iu$ in $\Afrei$ (also using the isomorphism of $\Iv\lhd\Aver$ and $\Toepl\otimes\K$ of Proposition \ref{DefIdealAverAtens}). We then compute the $K$-groups of $\Afrei$ using the following exact sequence:
\[0\to I\to \Afrei\to\Aver\to 0\]
In constrast, a natural first attempt to compute the $K$-theory of $\Afrei$ would try to make use of the following exact sequence:
\[0\to J\to \Afrei\to A_\theta\to 0 \]
Unfortunately it does not provide enough information in $K$-theory.

Throughout this section, $\theta\in\R$ is arbitrary (either rational or irrational).\\

We first recall the $K$-groups of $A_\theta$ and $\Aver$.

\begin{rem}
\begin{itemize}
\item[(a)]
The rotation algebra $A_\theta$ may be written as the crossed product  $\CS\rtimes_\theta\Z$ of $\CS$ with the automorphism $u\mapsto\lambda u$ on $\CS$, where $\lambda=e^{2\pi i\theta}$. We apply the Pimsner-Voiculescu sequence (\cite{PV80}) to obtain $K_0(A_\theta)=\Z^2$, generated by $1$ and a Rieffel projection (\cite{R81}), and $K_1(A_\theta)=\Z^2$, generated by the classes of the unitaries $u$ and $v$. (If $\theta=0$, then $K_0(A_\theta)=\Z^2$ is  generated by 1 and the Bott projection, cf. also \cite{Y86} for rational parameters $\theta$.)
\item[(b)]
 The $K$-groups of $\Aver$ are $K_0(\Aver)=\Z$, generated by the unit $1$, and $K_1(\Aver)=\Z$, generated by the class of the unitary $u$.
\end{itemize}
\end{rem}
\Proof (b) We compute the $K$-groups of $\Aver$ in exactly the same way as in the case of the rotation algebra. 
We use the Pimsner-Voiculescu sequence with respect to the automorphism on $\Toepl$ mapping $v\mapsto\lambda v$ on $\Toepl$. It induces the identity on the level of $K$-theory, hence the 6-term exact sequence splits into two parts:
\begin{align*}
 &&K_0(\Toepl) &&\starrow{0} &&K_0(\Toepl) &&\longrightarrow &&K_0(\Aver)\\
 &&\uparrow &&                               &&         &&    &&\downarrow\\
 &&K_1(\Aver)&&\longleftarrow &&K_1(\Toepl) &&\stackrel{0}{\longleftarrow} &&K_1(\Toepl) 
\end{align*}
Using $K_0(\Toepl)=\Z$ and $K_1(\Toepl)=0$ we get that $K_0(\Aver)=K_0(\Toepl)=\Z$ and $K_1(\Aver)=K_0(\Toepl)=\Z$. The generator of $K_0(\Toepl)=\Z$ is the unit $1$ and thus also of $K_0(\Aver)=\Z$. 

For the generator of $K_1(\Aver)=\Z$, we take a quick look at Pimsner and Voiculescu's proof of the 6-term exact sequence for $A\rtimes_\alpha\Z$, where $A$ is a unital $C^*$-algebra and $\alpha$ is an automorphism on $A$. They define $\Toepl(A,\alpha)$ to be the $C^*$-subalgebra of $(A\rtimes_\alpha \Z)\otimes \Toepl$ generated by all elements $a\otimes1$ for $a\in A$ and by $u\otimes v$, where $u$ is the adjoint unitary in $A\rtimes_\alpha\Z$. Then $A\otimes \K$ is isomorphic to the ideal in $\Toepl(A,\alpha)$ generated by $1\otimes(1-vv^*)$, mapping $a\otimes x_0$ to $a\otimes(1-vv^*)$ ($x_0$ being a minimal projection in $\K$). In consequence, the following sequence is exact:
\[0\to A\otimes\K\to\Toepl(A,\alpha)\to A\rtimes_\alpha\Z\to0\]
Now, Pimsner and Voiculescu show that $K_i(\Toepl(A,\alpha))$ is isomorphic to $K_i(A)$ and they conclude their picture of the 6-term exact sequence. The isometry $u\otimes v\in \Toepl(A,\alpha)$ is mapped to the unitary $u\in A\rtimes_\alpha\Z$.
Thus, the connecting map $\delta$ from $K_1(A\rtimes_\alpha\Z)$ to $K_0(A\otimes\K)$ maps the class of $u$ to the class of the defect projection $1\otimes(1-vv^*)$ of $u\otimes v$, which is isomorphic to $1\otimes x_0\in A\otimes\K$. This in turn is mapped to the unit of $A$ under the isomorphism of $K_0(A\otimes\K)$ and $K_0(A)$.
We conclude that the class of the unitary $u$ in $\Aver$ is mapped to the generator $1$ of $K_0(\Toepl)$, hence it must generate $K_1(\Aver)=\Z$.
\qqed

The $K$-groups of $\Atens$ are also obtained by a simple investigation of a 6-term exact sequence. 

\begin{thm}\label{KTheorieAtens}
 The $K$-groups of $\Atens$ are $K_0(\Atens)=\Z$, generated by the unit $1$ and $K_1(\Atens)=0$.
\end{thm}
\Proof
Consider the following exact sequence:
\[0\to\Iu\to\Atens\to\Aver\to0\]
By Proposition \ref{DefIdealAverAtens}, we know that the ideal $\Iu\lhd\Atens$ is isomorphic to $\Toepl\otimes\K$ via $u^iv(1-uu^*)\mapsto v\otimes x_i$ for $i\in\N_0$. 
The $K_0$-group of $\Toepl\otimes\K$ is $\Z$, generated by the class of $1\otimes x_0$, thus $K_0(\Iu)=\Z$ is generated by the class of the defect projection $1-uu^*$. But this is mapped to zero under  the map $K_0(\Iu) \to K_0(\Atens)$.
Indeed, the  class $[1-uu^*]\in K_0(\langle1-uu^*\rangle)$ is mapped to $[1-uu^*]\in K_0(\Atens)$. But  $[uu^*]=[u^*u]=1$ in $K_0(\Atens)$ which yields $1=[1-uu^*]+[uu^*]=[1-uu^*]+1$. Thus, the generator of $K_0(\Iu)$ is mapped to zero in $K_0(\Atens)$ and the map $K_0(\Iu)\to K_0(\Atens)$ is zero.

Furthermore $K_1(\Iu)=K_1(\Toepl)=0$, hence the 6-term exact sequence in $K$-theory corresponding to the above short exact sequence falls into the parts \linebreak $K_0(\Atens)\cong K_0(\Aver)=\Z$ (generated by the unit) and:
\[0\to K_1(\Atens)\to K_1(\Aver)\to K_0(\langle1-uu^*\rangle)\to 0\]
Since $K_1(\Aver)=\Z$ and $K_0(\langle 1-uu^*\rangle)=\Z$, we end up with $K_1(\Atens)=0$.
\qqed

A priori, the situation in the case of $\Afrei$ is less clear. The defect ideal \linebreak $\Iu\lhd\Afrei$ is more complicated than in the case of $\Atens$ and thus the $K$-groups cannot be computed in the same straightforward manner. But due to the embedding $\iota$ of $\Afrei$ into $\Asort$, we  are able to compute the $K$-groups of the ideal $J=\Iuv$ in $\Afrei$. Recall that $J$ is isomorphic to $\J$ (by Theorem \ref{JGestalt}), therefore we first consider the $C^*$-algebra $\Aup$.

In 1982, Cuntz computed the $K$-groups of the amalgamated free product of $C^*$-algebras (\cite{C82}). For the unital free product, we have the following result.

\begin{lem}[\cite{C82}]
 Let $A$ and $B$ be unital $C^*$-algebras and let $\psi_1:A\to \C$ and $\psi_2:B\to\C$ be unital homomorphisms. Then $K_0(A*_\C B)$ is the quotient of $K_0(A)\oplus K_0(B)$ by the subgroup generated by $([1_A],-[1_B])$.

Secondly, $K_1(A*_\C B)=K_1(A)\oplus K_1(B)$.
\end{lem}

Recall, that we view $\Aup$ as the universal $C^*$-algebra generated by a unitary $\bar u$ (corresponding to the unitary of $\CS$) and a projection $\bar p$ (corresponding to the vector $\left(\begin{array}{c} 1 \\ 0 \end{array}\right)\in\C^2$, cf. Lemma \ref{AupGestalt}). We apply the preceding lemma to $\Aup$.

\begin{prop}\label{KTheorieAup}
 The $K$-groups of $\Aup$ are $K_0(\Aup)=\Z^2$ generated by $[\bar p]$ and $[1-\bar p]$, and $K_1(\Aup)=\Z$ generated by $[\bar u]$.
\end{prop}
\Proof There are unital homomorphisms $\psi_1:\CS\to\C$, mapping $u\mapsto 1$, and $\psi_2:\C^2\to\C$ the projection onto the first component. By Cuntz' result, $K_0(\Aup)$ is the quotient of $K_0(\CS)\oplus K_0(\C^2)=\Z\oplus\Z^2$ by the subgroup generated by $\left(1,-\left(\begin{array}{c} 1 \\ 1 \end{array}\right)\right)$. But  $\Z\oplus\Z^2$ is generated by exactly this vector together with $\left(0,\left(\begin{array}{c} 1 \\ 0 \end{array}\right)\right)$ and $\left(0,\left(\begin{array}{c} 0 \\ 1 \end{array}\right)\right)$. So,  we end up with $K_0(\Aup)=\Z^2$, generated by $\left(\begin{array}{c} 1 \\ 0 \end{array}\right)$ and $\left(\begin{array}{c} 0 \\ 1 \end{array}\right)$ respectively $[\bar p]$ and $[1-\bar p]$.

Furthermore, $K_1(\Aup)=K_1(\CS)\oplus K_1(\C^2)=K_1(\CS)=\Z$. The generator is the class of  $u\in\CS$ respectively $\bar u\in\Aup$.\qqed

This yields the $K$-groups of $J$ in $\Afrei$.

\begin{cor}\label{KTheorieJ}
 The $K$-groups of the ideal $J=\Iuv$ in $\Afrei$ are $K_0(J)=\Z^2$ generated by the classes $[1-vv^*]$ and $[1-uu^*]$, and $K_1(J)=\Z$ generated by  $[u(1-vv^*)+(1-uu^*)v^*+uvv^*u^*]=[\bar u + (1-\bar u\bar u^*)]$. Here $\bar u$ is defined by $\bar u=u(1-vv^*)+(1-uu^*)v^*\in\Afrei$ as in Proposition  \ref{AfreiDecomp}.
\end{cor}
\Proof By Proposition \ref{KTheorieAup} we know that $K_0(\Aup)=\Z^2$ is generated by $[\bar p]$ and $[1-\bar p]$, thus $K_0(\J)=\Z^2$ is generated by $[\bar p\otimes x_0]$ and $[(1-\bar p)\otimes x_0]$. By the isomorphism of Theorem \ref{JGestalt} we infer that $\bar p\otimes x_0\in\J$ corresponds to $\bar p=v(1-uu^*)v^*\in J\subset \Afrei$, whereas $(1-\bar p)\otimes x_0$ corresponds to $p-\bar p=1-vv^*\in J$ (recall that $p$ is the unit of $C^*(\bar u, \bar p)\subset J$, see Remark \ref{WeiteresAfrei}(c)).
Because $v(1-uu^*)v^*$ is Murray von Neumann equivalent to $1-uu^*$ via $v(1-uu^*)$, we end up with $[1-vv^*]$ and $[1-uu^*]$ as generators of $K_0(J)=\Z^2$.

Now, $K_1(\J)=\Z$ is generated by $[\bar u\otimes x_0 + (1_{\widetilde{\J}}-1\otimes x_0)]$, where $\widetilde{\J}$ is the unitalization of $\J$. This is, since for any unital $C^*$-algebra $A$, the embedding $A\inj A\otimes \K$, $a\mapsto a\otimes x_0$ lifts to an isomorphism in $K$-theory.
Furthermore, let $A$ and $B$ be unital $C^*$-algebras, let \linebreak $\phi:A\to B$ be a homomorphism and $u\in M_n(A)$ a unitary. Then the class \linebreak $[u]\in K_1(A)$ is mapped under $\phi_*$ to $[\phi(u) + (1_n - \phi(1_n))]\in K_1(B)$ (see for instance \cite[exercise 8.5]{RLL00}). This is the canonical way of turning $K_1$ into a functor. It is consistent with the embeddings $U_nA\inj U_{n+1}A$ of unitary matrices $u\mapsto  
 \begin{pmatrix} u& \\ & 1 \end{pmatrix}$ by ``filling up with units''.

By Theorem \ref{JGestalt}, $\bar u\otimes x_0 + (1_{\widetilde{\J}}-1\otimes x_0)\in\J$ is isomorphic to $\bar u + (1-\bar u\bar u^*)\in\Afrei$, where $\bar u=u(1-vv^*)+(1-uu^*)v^*\in\Afrei$.
 Note, that $1-\bar u\bar u^*=uvv^*u^*$.\qqed

Next, we  compute the $K$-groups of the ideal $I:=\Iu$ in $\Afrei$. This will provide enough information to obtain the $K$-groups of $\Afrei$ from the following exact sequence:
\[0\to I\to \Afrei\to\Aver\to 0\]

\begin{prop}\label{KTheoryIAfrei}
 The $K$-groups of $I=\langle1-uu^*\rangle\lhd\Afrei$ are $K_0(I)=\Z$ generated by $[1-uu^*]$ and $K_1(I)=0$.
\end{prop}

\Proof Consider the following exact sequences of ideals in $\Afrei$ respectively $\Aver$ (recall the definitions $I=\langle 1-uu^*\rangle$ and $J=\langle 1-uu^*, 1-vv^*\rangle$):
\begin{align*}
 0&&\to &&I &&\to &&\Afrei &&\to &&\Aver &&\to &&0 \\
    &&  &&=  &&  &&\bigtriangledown && &&\bigtriangledown \\
 0&&\to &&I &&\to &&J &&\to &&\Iv &&\to &&0 
\end{align*}
As $I':=\Iv\lhd\Aver$ is isomorphic to $\CS\otimes\K$ by Proposition \ref{DefIdealAverAtens}, we know that $K_0(I')=\Z$ is generated by the class of $1-vv^*$, whereas $K_1(I')=\Z$ is generated by the class of $u(1-vv^*)+vv^*$. The latter uses the correspondence of $u\otimes x_0 + (1_{\widetilde{\CS\otimes\K}}-1\otimes x_0)\in\CS\otimes\K$ with $u(1-vv^*) + (1-(1-vv^*))\in\Aver$.
The second above sequence yields the following 6-term exact sequence in $K$-theory:
\begin{align*}
 &&K_0(I) &&\longrightarrow &&K_0(J) &&\longrightarrow &&K_0(I')\\
 && \uparrow &&   &&         &&    &&\downarrow\\
 &&K_1(I')&&\longleftarrow &&K_1(J) &&\longleftarrow &&K_1(I) 
\end{align*}
The map $K_0(J)=\Z^2\to K_0(I')=\Z$ maps one generator $[1-vv^*]$ of $K_0(J)$ to the generator $[1-vv^*]$ of $K_0(I')$. Thus it is surjective and hence the connecting map $K_0(I')\to K_1(I)$ is zero.
Secondly, the map $K_1(J)=\Z\to K_1(I')=\Z$ is an isomorphism, as the generator $[u(1-vv^*)+(1-uu^*)v^* + uvv^*u^*]\in K_1(J)$ is mapped to the generator $[u(1-vv^*)+vv^*]\in K_1(I')$. (Note that the element $u\in \Aver$ is a unitary and $uvv^*u^*=vv^*$ in $\Aver$.) 

 We conclude $K_1(I)=0$ and from
\[0\to K_0(I) \to K_0(J)=\Z^2 \to K_0(I')=\Z \to 0\]
we obtain $K_0(I)=\Z$. The generator is $[1-uu^*]$, which is mapped to the second generator $[1-uu^*]$ of $K_0(J)=\Z^2$.
\qqed

The result on $I\lhd\Afrei$ may also be read in terms of the ideal $\langle\bar p\rangle$ in $\Aup$.

\begin{rem}
 The ideal $I=\Iu\lhd\Afrei$ is isomorphic to the ideal $\langle\bar p\rangle\otimes\K$ in $\J$, mapping $v(1-uu^*)v^*$ to $\bar p\otimes x_0$.
Similarly, $\Iv\lhd\Afrei$ is isomorphic to $\langle1-\bar p\rangle\otimes\K\lhd\J$, where $1-vv^*$ is mapped to $(1-\bar p)\otimes x_0$.
 The $K$-groups of $\langle\bar p\rangle\lhd\Aup$ are $K_0(\langle\bar p\rangle)=\Z$, generated by $\bar p$ and $K_1(\langle\bar p\rangle)=0$.
Similarly, $K_0(\langle1-\bar p\rangle)=\Z$, generated by $1-\bar p$ and $K_1(\langle1-\bar p\rangle)=0$.
\end{rem}
\Proof Under the isomorphism of $J$ and $\J$ of Theorem \ref{JGestalt}, $(uv)^i\bar p=(uv)^iv(1-uu^*)v^*\in J$ is mapped to $\bar p\otimes x_i\in\J$ for $i\in \N_0$. Thus the image of $I=\langle v(1-uu^*)v^*\rangle$ in $\J$ is exactly $\langle\bar p\rangle\otimes\K$ (cf. also Proposition \ref{propKonkreteKorrespondenzen}). 

The ideal $\langle\bar p\rangle\lhd\Aup$ is the kernel of the map $\Aup\to\CS$, which maps $\bar u\mapsto u$ and $\bar p\mapsto 0$. Furthermore, the map $\CS\to\Aup$, mapping $u\mapsto\bar u$ is a split. Hence,  we have two exact sequences in $K$-theory, for $i=0,1$:
\[0\to K_i(\langle\bar p\rangle) \to K_i(\Aup)\to K_i(\CS)\to 0\]
As $K_i(\CS)=\Z$, these sequences are split exact which yields the result using Proposition \ref{KTheorieAup}.
\qqed

We are now able to compute the $K$-groups of $\Afrei$.

\begin{thm}\label{KTheoryAfrei}
 The $K$-groups of $\Afrei$ are $K_0(\Afrei)=\Z$, generated by the class of the unit $1$ and $K_1(\Afrei)=0$.
\end{thm}
\Proof Consider the following exact sequence, where $I=\Iu\lhd\Afrei$:
\[0\to I\to \Afrei\to\Aver\to0\]
 Under the map $K_0(I) \to K_0(\Afrei)$, the class of the projection $1-uu^*$ is mapped to zero, hence the map is zero.
 The proof for this is word-by-word the same as in the case of Theorem \ref{KTheorieAtens}.
Hence the 6-term exact sequence in $K$-theory corresponding to the above short exact sequence falls into the part $K_0(\Afrei)\cong K_0(\Aver)=\Z$,  where $[1]\in K_0(\Afrei)$ is mapped to $[1]\in K_0(\Aver)$, and the sequence
\[0\to K_1(\Afrei)\to K_1(\Aver)=\Z \to K_0(I)=\Z\to 0\]
from which we deduce $K_1(\Afrei)=0$.\qqed

\begin{rem}\label{MurphyKTheory}
Murphy (\cite{M94}) considered $C^*$-algebras associated to unital semigroups endowed with a 2-cocycle. He gave the example of $C^*_\theta(\N^2)$ in the case of the semigroup $\N^2$ and a 2-cocycle constructed out of a single complex scalar $\lambda$ of absolute value one. This is exactly the $C^*$-algebra $\Afrei$. In the introduction to his article, he mentions that the $K$-theory of this $C^*$-algebra was unknown, even in the trivial case of $C^*(\N^2)$, which is $\Afrei$ with $\theta=0$. According to Murphy, the knowledge of this $K$-theory would help in the theory of generalized Toeplitz operators (see \cite{M94} or \cite{M96} for references on this).
\end{rem}

\section*{Appendix}

\noindent
This appendix is supposed to give a brief overview on universal $C^*$-algebras. It includes some examples and some constructions of $C^*$-algebras viewed as universal $C^*$-algebras (cf. \cite[Sect. II.8.3]{B06} and \cite{C93}). \\

\zw{1}{The Construction of universal $C^*$-algebras} 
Let $G=\{x_i \, | \, i\in I\}$ be a set of \emph{generators} and denote by $P(G)$ the involutive algebra of non-commutative polynomials in $G\cup G^*$. Let $R\subset P(G)$ be a set of \emph{relations} and put $J(R)\subset P(G)$ to be the involutive ideal in $P(G)$ generated by the relations $R$. 
Then, the quotient $\mathcal A(G,R)$ of $P(G)$ by $J(R)$ is the universal involutive algebra, generated by $G$ and the relations $R$.
For all $x\in\mathcal A(G,R)$ we define:
\[\norm{x}:=\textnormal{sup}\{p(x)\,|\, p \textnormal{ is a $C^*$-seminorm on } \mathcal A(G,R) \}\]
A $C^*$-seminorm is a seminorm $p$ such that $p(xy)\leq p(x)p(y)$ and $p(x^*x)=p(x)^2$. 
If this expression $\norm{x}$ is finite for all $x\in\mathcal A(G,R)$, then $\norm{\cdot}$ is a $C^*$-seminorm on $\mathcal A(G,R)$. We take the quotient of $\mathcal A(G,R)$ by the null space $\{x\in\mathcal A(G,R) \, | \, \norm{x}=0\}$  and complete it with respect to $\norm{\cdot}$. We obtain the \emph{universal $C^*$-algebra generated by} $G$ \emph{and the relations} $R$, which we  denote by $C^*(G\;|\;R)$.

The most important feature of universal $C^*$-algebras is their \emph{universal property}. This is, for any $C^*$-algebra $B$ with elements $y_i\in B$, $i\in I$, such that the relations $R$ are fulfilled, there exists a unique $^*$-homomorphism from $C^*(G\;|\;R)$ to $B$, mapping $x_i\mapsto y_i$ for all $i\in I$.\\

\zw{2}{The rotation algebra} 
A classic example is the \emph{rotation algebra} $A_\theta$, where $\theta\in\R$ is a parameter of rotation. If $\theta$ is irrational, then $A_\theta$ is called \emph{irrational rotation algebra}, and \emph{rational rotation algebra} otherwise. We can define it as the following universal $C^*$-algebra:
\[A_\theta=\univbig{u,v\textnormal{ unitaries}}{uv=e^{2\pi i \theta} vu}\]
The algebraic relations for a unitary $u$ are $u^*u=uu^*=1$. (Note that $1$ is also one of the generators of $A_\theta$.)\\

\zw{3}{The Toeplitz algebra}
The \emph{Toeplitz algebra} $\Toepl$ may be viewed as the universal $C^*$-algebra generated by an isometry, thus:
\[\Toepl = \univ{v}{v \textnormal{ is an isometry, i.e. } v^*v=1}\]
The Toeplitz algebra is isomorphic to the $C^*$-subalgebra of $\LN$ generated by the unilateral shift $S\in\LN$, which is given by $Se_n=e_{n+1}$ for all $n\in\N_0$. Here $(e_n)_{n\in\N_0}$ is an orthonormal basis of the Hilbert space $\ell^2(\N_0)$. This isomorphism is due to Coburn's theorem (\cite{Co67}).\\


\zw{4}{The algebra of complex valued matrices}
 The $C^*$-algebra $M_n(\C)$ of complex valued $n\times n$-matrices  may be written as a universal $C^*$-algebra in the following two ways.
\begin{itemize}
 \item[(a)] $M_n(\C)$ is the universal $C^*$-algebra generated by elements $e_{ij}$ (named \emph{matrix units}) for $i,j\in\{1,\ldots,n\}$ and the relations $e_{ij}e_{kl}=\delta_{jk}e_{il}$ and $e_{ij}^*=e_{ji}$ for all $i,j,k,l\in\{1,\ldots,n\}$.
 \item[(b)] $M_n(\C)$ is the universal $C^*$-algebra generated by elements $x_1,x_2,\ldots$,$x_n$ and the relations $x_i^*x_j=\delta_{ij}x_0$ for all $i,j\in\{1,\ldots,n\}$. 
\end{itemize}
The first picture of $M_n(\C)$ uses the correspondence of $e_{ij}$ with the matrix
 in $M_n(\C)$ where all entries are zero except a unit $1$ at the $i$-th row and $j$-th column.
This is an isomorphism since both $C^*$-algebras are $n^2$-dimensional.
The universal $C^*$-algebras of (a) and (b) are isomorphic via $e_{ij}\mapsto x_ix_j^*$. Thus, $x_i$ corresponds to the matrix unit $e_{i0}$.
Note, that the element $x_0$ is a projection with $x_ix_0=x_i$ and that the $x_i$ are partial isometries.\\

\zw{5}{The algebra of compact operators}
The algebra $\K$ of compact operators  on a separable, infinite-dimensional Hilbert space may be written as  a universal $C^*$-algebra in the following two ways.
\begin{itemize}
 \item[(a)] $\K$ is the universal $C^*$-algebra generated by elements $e_{ij}$, where $i,j\in \N_0$  and the relations $e_{ij}e_{kl}=\delta_{jk}e_{il}$ and $e_{ij}^*=e_{ji}$ for all $i,j,k,l\in\N_0$.
 \item[(b)] $\K$ is the universal $C^*$-algebra generated by elements $x_i, i\in\N_0$ and the relations $x_i^*x_j=\delta_{ij}x_0$ for all $i,j\in\N_0$. 
Note, that the element $x_0$ is a projection with $x_ix_0=x_i$  and that the $x_i$ are partial isometries for all $i\in\N_0$.
\end{itemize}
The elements $e_{ij}$ correspond to the rank-one-operators matching the $j$-th basis vector $e_j$ to $e_i$ (where $(e_i)_{i\in\N}$ is an orthonormal basis of the Hilbert space). Denote by $\phi$ the surjection from the universal $C^*$-algebra $A$ of item (a) to $\K$ obtained by the universal property. It is an isomorphism, because the $C^*$-subalgebra  $M_n:=\univ{e_{ij}}{i,j=0,\ldots n-1}\subset A$ is isomorphic to  $M_n(\C)$ (as $M_n(\C)$ is simple). Therefore, the restriction of $\phi$ to $M_n$ is injective, again by simplicity of $M_n\cong M_n(\C)$. As the union of all $M_n$ for $n\in\N_0$ is dense in $A$, the map $\phi$ is isometric on a dense subset of $A$, thus it is isometric on the whole of $A$. We conclude that $\phi$ is an isomorphism.

Recall the connection between $\K$, $\Toepl$ and $\CS$.
The elements $x'_i:=v^i(1-vv^*)\in \Toepl$ fulfill the relations of $\K$.  As $\K$ is simple, the map $x_i\mapsto x_i'$ is injective. Now, $\K$ is the closed linear span of all elements $x_i x_j^*$,  $i,j, \in \N_0$ (the  $x_ix_j^*$ are the only monomials in $\K$). Hence its image in $\Toepl$ is exactly  the ideal $\Iv$ which is the closed linear span of all elements $x_i'x_j'^*=v^i(1-vv^*)(v^*)^j$. We obtain the following short exact sequence:
\[0\to\K\to\Toepl\to\CS\to 0\]

\zw{6}{(Full) group $C^*$-algebras}
We can also formulate some standard constructions of $C^*$-algebras in terms of universal $C^*$-algebras.
Let $G$ be a discrete group. The \emph{(full) group $C^*$-algebra} $C^*(G)$ may be viewed as the following universal $C^*$-algebra:
\[C^*(G)=\univ{u_g \textnormal{ unitaries for } g\in G}{u_{gh}=u_gu_h, u_{g^{-1}}=u_g^* \textnormal{ for all } g,h\in G}\]
We observe that the unitaries $(u_g)_{g\in G}$ encode the structure of the group $G$.

If $G=\Z$, we get $C^*(\Z)=\univ{u}{u \textnormal{ is a unitary}}=\CS$ since all  unitaries $u_k$ may be obtained as a power of $u_1$ or $u_1^*$. 
If $G=\Z_2$, the group of integers $\Z$ modulo $2\Z$, then $C^*(\Z_2)$ is the universal $C^*$-algebra generated by a self-adjoint unitary $u$ such that $u^2=1$.   Thus the only monomials  are of the form $1$ and $u$, and $C^*(\Z_2)$ is isomorphic to $\C^2$ via  $u\mapsto \left(\begin{array}{c} 1 \\ -1 \end{array}\right)$.\\

\zw{7}{Crossed products}
 Let $A$ be a unital $C^*$-algebra, let $G$ be a discrete group, and let $\alpha: G\to \textnormal{Aut}(A)$ be a group homomorphism.
The \emph{crossed product} of $A$ with $G$ may be seen as the following universal $C^*$-algebra:
\[A\rtimes_\alpha G = \univbig{a\in A, u_g \textnormal{ unitaries for } g\in G}{\textnormal{the relations } (*) \textnormal{ are fulfilled}}\]
The elements $a\in A$ fulfill the relations of $A$ and the relations $(*)$ are given by:
\[(*) \qquad u_{gh}=u_gu_h, u_{g^{-1}}=u_g^*, u_gau_g^*=\alpha_g(a) \textnormal{ for all } a\in A \textnormal{ and } g,h\in G\]
We infer that $A\rtimes_\alpha G$ is the $C^*$-algebra $A$ with adjoint unitaries which implement the automorphisms $\alpha_g$, $g\in G$.
If $A$ is not unital, then $A\rtimes_\alpha G$ is defined as the ideal, generated by all elements $a\in A$ in the above universal $C^*$-algebra.

Consider now a special case. Let $A$ be a unital $C^*$-algebra and let $\alpha$ be an automorphism on $A$. This gives rise to an action of $\Z$ on $\textnormal{Aut}(A)$, by $n\mapsto \alpha^n$. Hence we can form the crossed product of $A$ with $\Z$ which is given by:
\[A\rtimes_\alpha \Z=\univbig{a\in A, u \textnormal{ unitary}}{uau^*=\alpha(a) \textnormal{ for all } a\in A}\]
Note that we abuse the notation, because the automorphism as well as the action of $\Z$ on $\textnormal{Aut}(A)$ are denoted by $\alpha$.
As an example, consider the automorphism $v\mapsto e^{2\pi i \theta} v$ of $\CS$ where now $v$ denotes the generating unitary of $\CS$. Then $A_\theta=\CS\rtimes_\theta\Z$.

Another special case is $A=\C$. We let a discrete group $G$ act trivially on $A$ and we obtain $\C\rtimes G=C^*(G)$.\\

\zw{8}{Tensor products}
Let $A$ be a unital, nuclear $C^*$-algebra and let $B$ be another unital $C^*$-algebra. The \emph{tensor product} of $A$ and $B$ may be written as:
\[A\otimes B=\univbig{a\in A, b\in B}{ab=ba \textnormal{ for all } a\in A, b\in B}\]
Thus $A\otimes B$ is the universal $C^*$-algebra generated by all elements $a\in A$ (with the relations of $A$) and all $b\in B$ (with the relations of $B$) such that any $a\in A$ commutes with any $b\in B$. The units of $A$ and $B$ are identified. If $A$ and $B$ are not necessarily unital, then $A\otimes B$ is the ideal generated by all products $ab$ in the above $C^*$-algebra.\\

\zw{9}{Free products}
 Let $A$ and $B$ be unital $C^*$-algebras. The \emph{(unital) free product} $A*_\C B$ of $A$ and $B$ is defined as  the universal $C^*$-algebra generated by all elements $a\in A$ (with the relations of $A$) and all $b\in B$ (with the relations of $B$) under the identification of the units of $A$ and $B$. There are no further relations between elements from $A$ and those from $B$.

The constructions of the free product of $C^*$-algebras and the one for groups fit together in the following way. Let $G_1$ and $G_2$ be discrete groups and denote by $G_1*G_2$ their free product (of groups). Then the full group $C^*$-algebra $C^*(G_1*G_2)$ is isomorphic to the unital free product (of $C^*$-algebras) of $C^*(G_1)$ and $C^*(G_2)$.
(Cf. Blackadar's book on $K$-theory \cite[chapter V, 10.11.11(f)]{B98} for instance.)

\section*{Acknowledgements}

\noindent
My special thanks go to Joachim Cuntz, who encouraged me to work on isometries with twisted commutation relations. The discussions with him were very helpful and I appreciate his ideas and suggestions a lot. I would also like to thank the research group in Functional Analysis and Operator Algebras in M\"unster, in particular Siegfried Echterhoff, Thomas Timmermann, and Christian Voigt.

\end{document}